\documentclass[11pt]{amsart}
\usepackage{amsmath, amscd}

\newtheorem{Main}{Theorem}
\numberwithin{equation}{section}

\theoremstyle{plain}

\newtheorem{Prop}{Proposition}[section]
\newtheorem{Thm}[Prop]{Theorem}		
\newtheorem{Cor}[Prop]{Corollary}
\newtheorem{Lem}[Prop]{Lemma}
\newtheorem{Claim}{Claim}

\theoremstyle{definition}

\newtheorem{Def}[Prop]{Definition}
\newtheorem{Rem}[Prop]{Remark}

\newtheorem{Assume}[Prop]{\P}

\newcommand{\CBbb}{\mathbb C}

\newcommand{\ansing}{Z^{an.}}

\newcommand{\ansingtilde}{\widetilde Z^{an.}}

\newcommand{\algsing}{Z^{alg.}}

\newcommand{\Gr}{\mathop{\rm Gr}\nolimits}
\newcommand{\Grhn}{\mathop{\rm Gr}^{hn}\nolimits}
\newcommand{\Grs}{\mathop{\rm Gr}^{s}\nolimits}
\newcommand{\Grhns}{\mathop{\rm Gr}^{hns}\nolimits}

\newcommand{\filt}{\mathop{{\mathbb F}}^{hn}\nolimits}
\newcommand{\sat}{\mathop{\rm Sat}\nolimits}
\newcommand{\sing}{\mathop{\rm Sing}\nolimits}
\newcommand{\tor}{\mathop{\rm Tor}\nolimits}

\newcommand{\YM}{\mathop{\rm YM}\nolimits}
\newcommand{\HYM}{\mathop{\rm HYM}\nolimits}

\newcommand{\End}{\mathop{\rm End}\nolimits}

\newcommand{\vol}{\mathop{\rm vol}\nolimits}

\newcommand{\rk}{\mathop{\rm rk}\nolimits}

\newcommand{\tr}{\mathop{\rm Tr}\nolimits}
\newcommand{\dbar}{\bar\partial}
\newcommand{\lra}{\longrightarrow}
\newcommand{\weakarrow}{\rightharpoonup}
\newcommand{\diag}{\mathop{\rm diag}\nolimits}

\newcommand{\G}{\mathfrak G}

\newcommand{\GC}{{\mathfrak G}^{\CBbb}}

\newcommand{\A}{\mathcal{A}}
\newcommand{\Aone}{\mathcal{A}^{1,1}}

\begin{document}

% Topmatter

\title[The Yang-Mills Flow on K\"ahler Surfaces]
	{Convergence Properties of the Yang-Mills \\
Flow on K\"ahler Surfaces}

\author[Daskalopoulos]{ Georgios D. Daskalopoulos}

\address{Department of Mathematics \\
		Brown University \\
		Providence,  RI  02912}

\thanks{G.D. supported in part by NSF grant DMS-0204191}

\email{daskal@math.brown.edu}

\author[Wentworth]{ Richard A. Wentworth}

\address{Department of Mathematics \\
   Johns Hopkins University \\
   Baltimore, MD 21218}

\thanks{R.W. supported in part by NSF grant DMS-0204496}

\email{wentworth@jhu.edu}

%\subjclass{Primary: 57M50 ; Secondary: 58E20, 53C24}
\date{November 1, 2003}

\begin{abstract} 
Let $E$ be a  hermitian complex vector bundle over a compact K\"ahler surface $X$ with  K\"ahler form $\omega$,  and let $D$ be an integrable unitary connection on $E$ defining a holomorphic structure $D^{\prime\prime}$ on $E$.
We prove that the Yang-Mills flow on $(X,\omega)$  with initial condition $D$  converges, in an appropriate sense which takes into account bubbling phenomena,  to the double dual  of the graded sheaf associated to the $\omega$-Harder-Narasimhan-Seshadri filtration of the holomorphic bundle $(E,D^{\prime\prime})$.  This generalizes to K\"ahler surfaces the known result on Riemann surfaces and proves, in this case, a conjecture of Bando and Siu.
\end{abstract}

%$\hbox{}$
%\centerline{\fbox{\bf *** Preliminary Version ***}}
%\vskip .75in

\dedicatory{Dedicated to Professor Karen K. Uhlenbeck, \\ on the occasion of her 60th birthday.}

% End Topmatter

\maketitle

\thispagestyle{empty}

%\newpage

\baselineskip=16pt
\setcounter{footnote}{0}

 %%%%%%%%%%%%%%%%%%%%%%%%%%%%%%%%%%%%%%%%%%%%%%%%%%%%%%%%%%%%%%%%%%%%%%%%%%%%%%%%%%%%%%%%%%%%%%%%%%%%%%%%%%%%%%%%%%
%%%%%%%%%%%%%%%%%%%%%%%%%%%%%%%%%%%%%%%%%%%%%%%%%%%%%%%%%%%%%%%%%%%%%%%%%%%%%%%%%%%%%%%%%%%%%%%%%%%%%%%%%%%%%%%%%%

\section{Introduction}

 %%%%%%%%%%%%%%%%%%%%%%%%%%%%%%%%%%%%%%%%%%%%%%%%%%%%%%%%%%%%%%%%%%%%%%%%%%%%%%%%%%%%%%%%%%%%%%%%%%%%%%%%%%%%%%%%%%
%%%%%%%%%%%%%%%%%%%%%%%%%%%%%%%%%%%%%%%%%%%%%%%%%%%%%%%%%%%%%%%%%%%%%%%%%%%%%%%%%%%%%%%%%%%%%%%%%%%%%%%%%%%%%%%%%%

The main concern of  this paper is the relationship between the Yang-Mills functional and  stability of holomorphic vector bundles on K\"ahler surfaces.   By analogy with finite dimensional symplectic geometry on the one hand, and geometric invariant theory on the other, minimizing solutions to the Yang-Mills equations can be viewed as  zeros of a moment map on an infinite dimensional symplectic manifold, and the appearance of the stability condition for the existence of such solutions may be regarded as a version of the Kempf-Ness theorem.  More generally, one might expect a correspondence between two stratifications: the stable-unstable manifolds defined by the gradient flow of the Yang-Mills functional, and the algebraic stratification coming from  the maximally destabilizing one parameter subgroups.

This point of view  originated in the work of Atiyah and Bott \cite{AB} and was developed further by Donaldson \cite{Do1,Do2} (see \cite{Ki} for a general treatment of symplectic geometry vs.\ geometric invariant theory  in finite dimensions).
Let us recall some of the key points.
Given a holomorphic structure $\dbar_E$ on a complex vector bundle $E$ of rank $R$ over a K\"ahler manifold $X$,  one can associate a filtration by holomorphic subsheaves, called the Harder-Narasimhan filtration,  whose successive quotients are semistable.  The topological type of the pieces in the  associated graded object is encoded into an $R$-tuple $\vec \mu=(\mu_1,\ldots,\mu_R)$ of rational numbers called the Harder-Narasimhan type of $(E,\dbar_E)$.  Atiyah and Bott used the Harder-Narasimhan type to define a stratification of the infinite dimensional  space $\A^{\prime\prime}$ of holomorphic structures on $E$.  The group  $\GC$ of complex automorphisms of $E$, or the complex gauge group, acts on $\A^{\prime\prime}$ in a manner that preserves the stratification.  The main result of \cite{AB} is that, when $X$ is a Riemann surface,  the stratification obtained from the Harder-Narasimhan  filtration is $\GC$-equivariantly perfect, and this leads to a recursive calculation  of the cohomology of the moduli spaces of stable bundles, in certain cases.

On the other hand, fixing a hermitian structure $H$ on $E$, one may identify $\A^{\prime\prime}$ with the space $\A_H$ of unitary connections on $E$ via the map which sends a unitary connection $D$ to its (0,1) part $D^{\prime\prime}=\dbar_E$ (in higher dimensions we require the integrability condition that the curvature $F_D$ be of type (1,1)).
The Yang-Mills functional, which associates to a connection the $L^2$-norm of its curvature, can be used as a Morse function on $\A_H$.  For a fixed holomorphic bundle $(E,\dbar_E)$, the Harder-Narasimhan type gives an absolute lower bound on the Yang-Mills number of any connection in the isomorphism class.  Up to a topological term, the Yang-Mills number is the same  as the Hermitian-Yang-Mills number, which is defined as the $L^2$ norm of the contraction $\Lambda F_D$ of the curvature with the K\"ahler form (see \eqref{E:chernclass}).   If $\vec\mu$ is the Harder-Narasimhan type of $(E,\dbar_E)$, then:
$$
\HYM(\vec \mu):= 2\pi \sum_{i=1}^R\mu_i^2\leq \HYM(D):=\int_X|\Lambda F_D|^2\, dvol
$$
for all $D$ such that $(E,D^{\prime\prime})$ is holomorphically isomorphic to $(E,\dbar_E)$ (see Cor.\ \ref{C:minimum} for a proof of this result for K\"ahler surfaces).  Atiyah and Bott conjectured that, on a Riemann surface,  the gradient flow of the Yang-Mills functional should converge at infinity, achieving the lower bound expressed above.  Moreover, the stable-unstable manifold stratification should coincide with the Harder-Narasimhan stratification.  

That this is indeed the case follows from the work of several authors.  First, Donaldson proved the long time existence of the $L^2$-gradient flow for the Yang-Mills functional on any K\"ahler manifold \cite{Do1, DoKr}.  For Riemann surfaces, the  asymptotic  convergence of the gradient flow and the equivalence of the two stratifications was established in \cite{D} (see also \cite{R} for a more analytic approach).  A key fact which makes the two dimensional case more tractable is that the Yang-Mills functional on Riemann surfaces satisfies the equivariant Palais-Smale Condition C.

In higher dimensions Condition C fails.  More seriously, the Yang-Mills flow can develop singularities in finite time.  Therefore, one cannot expect  a Morse theory in the classical sense.  For Riemannian four-manifolds, Taubes formulated an extended  Morse theory by attaching bundles with varying topologies in order to compensate for curvature concentration, or \emph{bubbling}.  In this approach one uses the strong gradient of the Yang-Mills functional associated to a complete Riemannian metric on the space $\A_H$, which exists for all time by the fundamental existence theorem for ODE's.  In this way Taubes established  the connectivity of the moduli space of self dual connections
in certain cases.  Moreover, he was able to calculate the stable homotopy groups  as predicted by the conjecture of Atiyah and Jones \cite{T}.

In the case of K\"ahler surfaces it is more natural to consider the $L^2$ rather than the strong gradient flow.  Long time existence, as mentioned above, is guaranteed.  Donaldson used this flow to prove the correspondence between anti-self-dual connections and stable bundles  (see \cite{Do1}, and more generally \cite{Do2} and \cite{UY} for holomorphic bundles in higher dimensions).   This is now known as the Hitchin-Kobayashi Correspondence, or the Donaldson-Uhlenbeck-Yau Theorem.
In \cite{BS},   Bando and Siu 
extended the correspondence to coherent analytic sheaves by considering singular hermitian metrics with controlled curvature.  They also conjectured that the  relationship between the Yang-Mills flow and the Harder-Narasimhan filtration  which holds on Riemann surfaces should analogously be true in higher dimensions.    We use the word ``analogous"  because even when considering vector bundles  the Harder-Narasimhan filtration in higher dimensions may be given only by subsheaves, and the associated graded objects may not be locally free.  Indeed, it is for this reason that the generalization of Bando-Siu naturally leads to a  conjecture on the behavior of the flow for vector bundles.

The purpose of this paper is to prove the conjecture of Bando and Siu for holomorphic bundles on K\"ahler surfaces.
To state the result precisely, let $\Grhns_\omega(E,\dbar_E)$ denote the Harder-Narasimhan-Seshadri filtration of $(E,\dbar_E)$ with respect to the K\"ahler form $\omega$, and let $\Grhns_\omega(E,\dbar_E)^{\ast\ast}$ be its double dual. 
 To clarify, we note here that the Harder-Narasimhan-Seshadri filtration is actually a double filtration which takes into account the possibility that the successive factors in the Harder-Narasimhan filtration may only be semistable as opposed to stable (see Prop.\ \ref{P:hnsfiltration}).  Thus, the individual factors in the associated graded object  are all stable.
 Since $X$ is now assumed to be a surface, the double dual is a vector bundle and carries a Yang-Mills connection which realizes the hermitian structure as a direct sum of Hermitian-Einstein metrics.
  Now for unstable bundles the flow may not converge in the usual sense; again, because of bubbling.  However, one can always extract subsequential \emph{Uhlenbeck limits} which are Yang-Mills connections on bundles with a possibly different topology than the original $E$.  The bundles are isometric, and the connections converge, away from a singular  set of codimension four\footnote{\, In this paper, convergence of connections will always be modulo real gauge equivalence.}.  In dimension four, the singular set is a finite collection of points.  For the precise definition, see Prop.\ \ref{P:convergence} below.
 Our result is that on a K\"ahler surface, the Uhlenbeck limits are independent of the subsequence and  are determined solely by the isomorphism class of the  initial holomorphic bundle $(E,\dbar_E)$.  More precisely: 

\begin{Main}[Main Theorem]  \label{T:main}
Let $X$ be a compact K\"ahler surface, $E\to X$ a hermitian vector bundle, and $D_0$ an integrable unitary connection on $E$ inducing a holomorphic structure $\dbar_E=D_0^{\prime\prime}$.  Let $D_\infty$ denote the Yang-Mills connection on $\Grhns_\omega(E,\dbar_E)^{\ast\ast}$ referred to above.
Let $D_{t}$ be the time $t$ solution to the Yang-Mills flow with initial condition $D_0$.  Then as $t\to \infty$, $D_t$ converges in the sense of Uhlenbeck to  $D_\infty$.
\end{Main}

We now give a sketch of the ideas involved in the  proof of the Main Theorem and explain the organization of the paper.  In Section \ref{S:prelim} we lay out the definitions of the Harder-Narasimhan-Seshadri filtration and its associated graded object.  We review the Yang-Mills flow and the notion of an Uhlenbeck limit.  We also discuss other Yang-Mills type functionals associated to invariant convex functions on the Lie algebra of the unitary group.  These are closely related to $L^p$ norms, and they will play an  important role not only in distinguishing  the various strata, but also because one actually cannot expect good $L^2$ behavior in the constructions that follow.

Much of the  difficulty in proving the Main Theorem arises from the fact that the Harder-Narasimhan filtration is not necessarily given by subbundles.  Since we have restricted our attention in this paper to surfaces, the individual factors in the filtration are themselves locally free, but the successive quotients may have point singularities.  These points are essentially the locus where one can expect bubbling to occur along the flow --  when the filtration is by subbundles one can show there is no bubbling  -- and they are therefore a fundamental aspect of the problem and not a mere technical annoyance. 

In Section \ref{S:blowup} we therefore analyze the degree to which the singularities in the filtration can be resolved by blowing up.  To be more precise, we are interested in comparing the Harder-Narasimhan filtrations of $E\to X$ and $\pi^\ast(E)\to \widehat X$, where $\pi: \widehat X\to X$ is a sequence of monoidal transformations, and the K\"ahler metric $\omega_\varepsilon$ on $\widehat X$ is an $\varepsilon$-perturbation by the components of the exceptional divisor of the pullback of the K\"ahler metric $\omega$ on $X$.  In Thm.\ \ref{T:hnresolution} we prove that under the assumption that the successive quotients of the Harder-Narasimhan filtration are stable,  there is a resolution $\widehat X$ such that for sufficiently small $\varepsilon$, the Harder-Narasimhan filtration of $\pi^\ast(E)$ is given by subbundles and its direct image by $\pi$  coincides with the Harder-Narasimhan filtration of $E$.  The situation is more complicated for semistable factors, and the resolution of the filtration by subbundles may not correspond to the Harder-Narasimhan filtration for any $\varepsilon>0$.  

Nevertheless,
this analysis is  sufficient for our purposes.  In particular, 
 we introduce the notion of an $L^p$-\emph{approximate critical hermitian structure}.  Roughly speaking, this is a smooth hermitian metric on a holomorphic bundle whose curvature in the direction of the K\"ahler form is close in the $L^p$ sense to a critical value determined by its Harder-Narasimhan type (see  Def.\ \ref{D:approximate}).  We prove Thm.\ \ref{T:approximate} which states that there exist $L^p$-approximate critical hermitian structures for all $1\leq p<\infty$.  This result is an $L^p$ version of a conjecture attributed to Kobayashi.  Interestingly, the method does not seem to extend to $p=\infty$.
 
 The first step in the proof of the Main Theorem is to determine the Harder-Narasimhan type of an Uhlenbeck limit.  Since the  Hermitian-Yang-Mills numbers are monotone along the flow one can show that for an initial condition which is a sufficiently close approximate critical hermitian structure,  the Uhlenbeck limit of a sequence along the  flow must have the correct Harder-Narasimhan type. Then  a length decreasing argument for the Yang-Mills flow, which closely resembles  Hartman's result for the harmonic map flow,  implies that \emph{any}  initial condition  must have Uhlenbeck limits of the correct type (see Thm.\ \ref{T:type}).

The second step in the proof of the Main Theorem is to show that the holomorphic structure on the Uhlenbeck limit coincides with the double dual of the associated graded sheaf of the Harder-Narasimhan-Seshadri filtration.  The approach here is necessarily completely different from that of \cite{D}.
The main idea is to generalize an argument of Donaldson who constructs limiting  holomorphic maps from the sequence of complex gauge transformations defined by the sequence of connections along the flow.  Instead of a map on the entire initial bundle, we show that maps can be formed for each of the pieces in the filtration separately.  The result then follows by an inductive argument.  

The proof that the limiting holomorphic structure is the correct one is largely independent of the details of the flow.  Indeed, we only use the fact that the Yang-Mills numbers of a descending sequence in a complex gauge orbit are absolutely minimizing.  Since this result is in some way disjoint from  Thm.\ \ref{T:main}, we formulate it separately (see Thm.\ \ref{T:holomorphic}).  Recall that sequences minimize Yang-Mills energy $\YM(D)$  if and only if they minimize Hermitian-Yang-Mills energy $\HYM(D)$:

\begin{Main}[Minimizing Sequences]
Let $X$ be a compact K\"ahler surface, $E\to X$ a hermitian vector bundle, and $D_0$ an integrable unitary connection on $E$ inducing a holomorphic structure $\dbar_E=D_0^{\prime\prime}$. Let $\vec\mu_0$ be the Harder-Narasimhan type of $(E, \dbar_E)$, and  let $D_\infty$ denote the Yang-Mills connection on $\Grhns_\omega(E,\dbar_E)^{\ast\ast}$.
Suppose $D_j$ is a sequence of smooth unitary connections in the complex gauge orbit of $D_0$ such that
 $\HYM(D_j)\to \HYM(\vec\mu_0)$ as $j\to \infty$.
Then there is a subsequence (also denoted $j$) and a finite set of points $\ansing\subset X$ such that
\begin{enumerate}
\item  $E$ and $\Grhns_\omega(E,\dbar_E)^{\ast\ast}$ are $L^p_{2,loc.}$-isometric on $X\setminus \ansing$ for all $p$;
\item $D_j\to D_\infty$ in $L^2_{loc.}$ away from $\ansing$.
\end{enumerate}
\end{Main}

\smallskip
\noindent
\emph{Acknowledgement.}   This paper is a substantially revised version of an earlier preprint.  We thank the referee for an exceptionally careful reading of that initial manuscript and for numerous and helpful comments.  The references \cite{Bu1,Bu2,Bu3} suggested by the referee, in particular, simplified some of our arguments and allowed us to remove the restriction to projective surfaces required in the original paper.

 \bigskip

%%%%%%%%%%%%%%%%%%%%%%%%%%%%%%%%%%%%%%%%%%%%%%%%%%%%%%%%%%%%%%%%%%%%%%%%%%%%%%%%%%%%%%%%%%%%%%%%%%%%%%%%%%%%%%%%%%
%%%%%%%%%%%%%%%%%%%%%%%%%%%%%%%%%%%%%%%%%%%%%%%%%%%%%%%%%%%%%%%%%%%%%%%%%%%%%%%%%%%%%%%%%%%%%%%%%%%%%%%%%%%%%%%%%%

\section{Preliminaries}       \label{S:prelim}

%%%%%%%%%%%%%%%%%%%%%%%%%%%%%%%%%%%%%%%%%%%%%%%%%%%%%%%%%%%%%%%%%%%%%%%%%%%%%%%%%%%%%%%%%%%%%%%%%%%%%%%%%%%%%%%%%
%%%%%%%%%%%%%%%%%%%%%%%%%%%%%%%%%%%%%%%%%%%%%%%%%%%%%%%%%%%%%%%%%%%%%%%%%%%%%%%%%%%%%%%%%%%%%%%%%%%%%%%%%%%%%%%%%%

%%%%%%%%%%%%%%%%%%%%%%%%%%%%%%%%%%%%%%%%%%%%%%%%%%%%%%%%%%%%%%%%%%%%%%%%%%%%%%%%%%%%%%%%%%%%%%%%%%%%%%%%%%%%%%%%%%

\subsection{Stability and the Harder-Narasimhan Filtration}  \label{S:hnfiltration}

Let $X$ be a complex surface.    The  singular set $\sing(E)$ of a coherent analytic torsion-free sheaf $E\to X$ is the closed subvariety where $E$ fails to be locally free.  Since $\dim_\CBbb X=2$, the singular set  of a torsion-free sheaf is a locally finite collection of points and reflexive sheaves are locally free  (cf.\ \cite[Cor.\  V.5.15 and V.5.20]{Ko}).
A subsheaf $S\subset E$ of a reflexive sheaf $E$ is  said to be saturated if the quotient $Q=E/S$ is torsion-free.  In general,
the saturation of a subsheaf $S$ in $E$, denoted $\sat_E(S)$, is the kernel of the sheaf map $E\to Q/\tor(Q)$, where $\tor(Q)$ is the torsion subsheaf of $Q$.
Note that $S$ is a subsheaf of $\sat_E(S)$ with a torsion quotient.  A saturated subsheaf of a reflexive sheaf is reflexive 
(cf.\
\cite[Prop.\ V.5.22]{Ko}).
We will also need the following result, whose proof is standard:

\begin{Lem} \label{L:sat}
Let $E$ be a torsion-free sheaf.  Suppose $S_1\subset S_2\subset E$ are subsheaves with $S_2/S_1$ a torsion sheaf.  Then $\sat_E(S_1)=\sat_E(S_2)$.
\end{Lem}

Now assume that $X$ is compact  with a  K\"ahler form $\omega$.  We will  assume the volume of $X$ with
respect to $\omega$ is normalized to be $\vol(X)=2\pi$. The $\omega$-slope $\mu(E)$  of a torsion-free sheaf
$E\to X$ is defined by:
\begin{equation} \label{E:slope}
\mu_\omega(E)=\frac{\deg_\omega (E)}{\rk (E)}=\frac{1}{\rk(E)}\int_X c_1(E)\wedge\omega\ .
\end{equation}
 
 We define $\mu_{max}(E)$ to be the maximal slope of a  subsheaf of $E$, and
$\mu_{min}(E)$ to be the minimal slope of a torsion-free quotient of $E$. 
A torsion-free sheaf $E\to X$ is $\omega$-stable (resp.\ $\omega$-semistable) if for all subsheaves $F\subset E$ with $0<\rk(F)<\rk(E)$,  $\mu_\omega(F) <
\mu_\omega(E)$  (resp.\ $\mu_\omega(F) \leq
\mu_\omega(E)$).
When the K\"ahler form is understood we shall sometimes  refer to $E$ simply as stable or
semistable, and we will also omit  subscripts  and write $\mu(E)$.  

\begin{Prop}[cf.\ \cite{Ko}, Thm.\ V.7.15]  \label{P:hnfiltration}
Let $E\to X$ be a torsion-free sheaf.  Then there is a filtration:
$
0=E_0\subset E_1\subset\cdots\subset E_\ell=E
$,
 called the  Harder-Narasimhan filtration of $E$ (abbr.\ 
HN filtration),
such that $Q_i=E_i/E_{i-1}$ is torsion-free and semistable.  Moreover, $\mu(Q_i)>\mu(Q_{i+1})$, and the associated graded object
$
\Grhn_{\omega}(E)=\oplus_{i=1}^\ell Q_i
$
is uniquely determined by the isomorphism class of $E$
  \end{Prop}

It will be convenient to denote the
subsheaf $E_i$ in the   HN filtration by $\filt_i(E)$, or by $\filt_{i,\omega}(E)$, when we wish to emphasize the role of the K\"ahler structure.
The collection of slopes $\mu(Q_i)$ is an important invariant of the isomorphism class of a torsion-free sheaf.  For a torsion-free sheaf $E$ of rank $R$  construct an $R$-tuple of numbers $\vec\mu(E)=(\mu_1,\ldots, \mu_R)$ from the HN filtration by setting:
 $\mu_i=\mu(Q_j)$, for $\rk(E_{j-1})+1\leq i\leq \rk(E_j)$.
  We call $\vec\mu(E)$ the 
Harder-Narasimhan type of $E$.
These invariants admit a natural partial
ordering which will be very relevant to this paper.  
 For a pair $\vec\mu$, $\vec\lambda$ of $R$-tuple's satisfying 
$\mu_1\geq\cdots\geq\mu_R$,  $\lambda_1\geq
\cdots\geq\lambda_R$,  and $\sum_{i=1}^R\mu_i=\sum_{i=1}^R\lambda_i$, we define:
\begin{equation} \label{E:slopeordering}
\vec\mu\leq\vec\lambda \qquad \iff\qquad \sum_{j\leq k}\mu_j\leq \sum_{j\leq k}\lambda_j\ ,\qquad\text{for all}\ k=1,\ldots, R\ .
\end{equation}
The importance of this ordering is that it defines a stratificaton of the space of holomorphic structures on a given complex vector bundle over a
Riemann surface.  See \cite[\S 7]{AB} for more details.  We will make use of the following  simple fact:
\begin{Lem} \label{L:shatz}
 Let $\vec\mu=(\mu_1,\ldots,\mu_R)$ and $\vec\lambda=(\lambda_1,\ldots,\lambda_R)$ be nonincreasing $R$-tuples as above.  Suppose there is a partition $0=R_0<R_1<\cdots < R_\ell=R$ such that $\mu_i=\mu_j$ for all pairs $i,j$ satisfying: 
 $
 R_{k-1}+1\leq i,j \leq R_k\ ,\ k=1,\ldots,\ell
 $.
   If 
$ \sum_{j\leq R_k}\mu_j\leq \sum_{j\leq R_k}\lambda_j$, for all $k=1,\ldots,\ell$, then $\vec\mu\leq \vec\lambda$.
 \end{Lem}

Several technical properties of the HN filtration will also play a role in this paper. We again omit the proofs.

\begin{Prop} \label{P:general}  
\begin{enumerate}
\item
Let $\widetilde E\to X$ be torsion-free and $E\subset\widetilde E$ with $T=\widetilde E/E$ a torsion sheaf supported at points.  Then 
$
\filt_i(E)=\ker(\filt_i(\widetilde E)\to T)
$, and
$
 \filt_i(\widetilde E)=\sat_{\widetilde E}(\filt_i(E))
$.
\item
Let $E\to X$ be a torsion-free sheaf.  Let $E_1=\filt_1(E)$, and $Q_1=E/E_1$.  Then
$$
\filt_{i+1}(E)=\ker( E\lra Q_1/\filt_i(Q_1))\ .
$$
In particular, $\filt_{i+1}(E)/E_1=\filt_i(Q_1)$.
\item
Consider an exact sequence:
$
0\to S\to E\to Q\to 0
$
of torsion-free sheaves on $X$ with $
\mu_{min}(S)>\mu_{max}(Q)$.
Then the Harder-Narasimhan filtration of $E$ is given by:
$$
0=E_0\subset \filt_1(S)\subset\cdots\subset\filt_k(S)=
S\subset\filt_{k+1}(E)\subset\cdots\subset\filt_\ell(E)=E
\ , 
$$
where $
\filt_i(E)=\filt_i(S)$,  $i\leq k$, and
$
\filt_{k+i}(E)=\ker(E\to Q/\filt_i(Q))$ for $ i=0,1,\ldots, \ell-k$.
In particular,
$
\Grhn(E)\simeq \Grhn(S)\oplus\Grhn(Q)
$.
\end{enumerate}
\end{Prop}

We point out the analogous filtrations for semistable sheaves:

\begin{Prop}[cf.\ \cite{Ko}, Thm.\ V.7.18]  \label{P:sfiltration}
Let $Q\to X$ be a semistable torsion-free sheaf.  Then there is a filtration
$
0=F_0\subset F_1\subset\cdots\subset F_\ell=Q
$,
called a Seshadri filtration of $E$,
such that $F_i/F_{i-1}$ is stable and torsion-free.  Moreover, $\mu(F_i/F_{i-1})=\mu(Q)$ for each $i$.  The associated graded object
$
\Grs_{\omega}(Q)=\oplus_{i=1}^\ell F_i/F_{i-1}
$,
is uniquely determined by the isomorphism class of $Q$.
\end{Prop}
 
 Finally, the double filtration whose associated graded sheaf appears in the statements of the main results in the Introduction is obtained by combining Prop.'s \ref{P:hnfiltration} and \ref{P:sfiltration}:
 
\begin{Prop}  \label{P:hnsfiltration}
Let $E\to X$ be a torsion-free sheaf.  Then there is a double filtration $\{E_{i,j}\}$,
 called a Harder-Narasimhan-Seshadri filtration of $E$ (abbr.\
HNS-filtration),  with the following properties:
if $\{E_i\}_{i=1}^\ell$ is the  HN filtration of $E$, then 
$
E_{i-1}=E_{i,0}\subset E_{i,1}\subset\cdots\subset E_{i,\ell_i}=E_i
$,
and the successive quotients $Q_{i,j}=E_{i,j}/E_{i,j-1}$ are stable torsion-free sheaves.  Moreover,
$
\mu(Q_{i,j})=\mu(Q_{i,j+1})
 $ and $\mu(Q_{i,j})>\mu(Q_{i+1,j})$.
The associated graded object:
$$
\Grhns_{\omega}(E)=\bigoplus_{i=1}^\ell\bigoplus_{j=1}^{\ell_i} Q_{i,j}
$$
is uniquely determined by the isomorphism class of $E$.
 \end{Prop}

%%%%%%%%%%%%%%%%%%%%%%%%%%%%%%%%%%%%%%%%%%%%%%%%%%%%%%%%%%%%%%%%%%%%%%%%%%%%%%%%%%%%%%%%%%%%%%%%%%%%%%%%%%%%%%%%%%

\subsection{Yang-Mills Connections and Uhlenbeck Limits}  \label{S:ym}

 Given a smooth complex vector bundle $E\to X$  of rank $R$,  let $\Omega^{p,q}(E)$ denote the space of smooth $(p,q)$ forms with values in $E$.
We will regard a holomorphic structure on $E$ as given by a $\dbar$ operator $\dbar_E:\Omega^{p,q}(E)\to \Omega^{p,q+1}(E)$ satisfying the integrability
condition  $\dbar_E\circ\dbar_E=0$.   We will sometimes denote the holomorphic structure on $E$ explicitly
by $(E,\dbar_E)$.  When this structure is understood, we will confuse the notation for the holomorphic
bundle and the sheaf of holomorphic sections by $E$, as we have done in the previous section.

Now suppose we are given a smooth hermitian metric $H$ on $(E,\dbar_E)$.  Then there is a uniquely determined $H$-unitary connection $D$ on
$E$ satisfying $D^{\prime\prime}=\dbar_E$, where $D^{\prime\prime}$ denotes the $(0,1)$ part of $D$ ($D'$ will denote the $(1,0)$ part).  We will sometimes
denote this connection by
$D=(\dbar_E, H)$.  Conversely, given a unitary connection
$D$ on
$E$ whose curvature $F_D=D\circ D$ is of type $(1,1)$ (i.e.\ $F_D^{0,2}=0$), then $D^{\prime\prime}=\dbar_E$ defines a
holomorphic structure on $E$, and $D=(D^{\prime\prime},H)$.

For a fixed hermitian metric,
let $\A_H$ denote the space of $H$-unitary connections $D$ on $E$, and let $\Aone_H$ denote those satisfying $F_D^{0,2}=0$.  The discussion above
gives an identification of $\Aone_H$ with the space
$\A^{\prime\prime}$ of integrable $\dbar$-operators, or holomorphic structures, on $E$. 
 We denote by $\G$ the space of unitary gauge transformations acting on $\A_H$ by pulling back.   Via the identification
$\Aone_H\simeq
\A^{\prime\prime}$,  it is clear that we have an action on $\Aone_H$ by
$\GC$, the complexification of $\G$.  We call $\GC$ the complex gauge group.  Notice that $\GC$ also acts on the space of hermitian metrics
on $E$, where $g(H)$ is defined by: $g(H)(s_1,s_2)=H(gs_1,gs_2)$.   
 
Since many norms will be used in this paper, let us emphasize the following:  if ${\bf a}$ is a hermitian or skew-hermitian endomorphism on an $R$-dimensional hermitian vector space with eigenvalues $\{\lambda_1,\ldots,\lambda_R\}$, we set:
 \begin{equation} \label{E:hermitiannorm}
 |{\bf a}|=\left\{ \sum_{i=1}^R |\lambda_i|^2\right\}^{1/2}\ .
 \end{equation}
 
 \noindent  For a hermitian vector bundle $E$, let ${\mathfrak u}(E)$ denote the subbundle of $\End E$ consisting of skew-hermitian endomorphisms.  If $\bf a$ is a section of ${\mathfrak u}(E)$, then $|{\bf a}|$ will denote the pointwise norm defined by (\ref{E:hermitiannorm}).

Given a K\"ahler metric $\omega$ on $X$,  the  Yang-Mills Functional (abbr.\ \emph{YM Functional}) is defined by
 $\YM(D)=\Vert F_D\Vert^2_{L^2(\omega)}$, and the  Hermitian-Yang-Mills Functional (abbr.\ \emph{HYM Functional}) is defined by $\HYM(D)=\Vert\Lambda_\omega F_D\Vert^2_{L^2(\omega)}$.
Here, $\Lambda_\omega$ denotes contraction with the K\"ahler form, and $\Lambda_\omega F_D$ is called the  Hermitian-Einstein tensor associated to $D$.
 Since for any $D\in\Aone_H$ we have (cf.\
\cite[IV.3.29]{Ko}):
\begin{equation} \label{E:chernclass}
\YM(D)=\HYM(D)+4\pi^2(2c_2(E)-c_1^2(E))\ ,
\end{equation}
 the YM and HYM functionals have the
same critical points on $\Aone_H$; namely, the  Yang-Mills Connections $D^\ast F_D=0$.  We also note  the K\"ahler identities:
\begin{equation} \label{E:kahleridentities}
D^\ast F_D =\sqrt{-1}\, (D'-D'')\Lambda_\omega F_D\ .
\end{equation}

Given a holomorphic bundle $E\to (X,\omega)$, a hermitian metric $H$ is called a Hermitian-Einstein Metric if there is a constant $\mu$ such that
$
\sqrt{-1}\, \Lambda_\omega F_{(\dbar_E,H)} = \mu\, {\bf I}_E
$, where ${\bf I}_E$ denotes the identity endomorphism of $E$.
If  $X$ is compact, then because of the normalization $\vol(X)=2\pi$  it is necessarily the case that $\mu=\mu_\omega(E)$ (see (\ref{E:slope})).
The celebrated theorem of Donaldson-Uhlenbeck-Yau  relates stability to the existence of a Hermitian-Einstein metric (cf.
\cite{Do1,Do2,UY}):

\begin{Thm} \label{T:DUY}
A holomorphic vector bundle $E$ on a compact K\"ahler manifold $(X,\omega)$ admits a Hermitian-Einstein metric if and only if it is
holomorphically split into a direct sum of $\omega$-stable bundles, all with slope $=\mu(E)$.
\end{Thm}

The following is standard  (cf.\ \cite[\S 4]{AB} or \cite[IV.
\S 3]{Ko}):

\begin{Prop} \label{P:split}
Let $D\in\Aone_H$ be a YM connection on a hermitian vector bundle over a K\"ahler manifold $X$.  Then there is an orthogonal splitting $(E,D)=\oplus_{i=1}^\ell (Q_i, D_i)$, where  
$\sqrt{-1}\,
\Lambda_\omega F_{D_i}=\mu_i {\bf I}_{Q_i}$, for constants $\mu_i$.  In case $X$ is compact, $\mu_i=\mu(Q_i)$,  and the
critical values of the YM functional on $\Aone_H$ are discrete.
\end{Prop}

 For the proof of the following
version of Uhlenbeck compactness, see \cite{U1} (and also, \cite[Thm.\ 5.2]{UY}):
\begin{Prop} \label{P:compactness}
Let $X$ be a compact K\"ahler surface  and $E\to X$ a complex vector bundle with hermitian metric
$H$.  Assume
$D_j$ is a sequence of integrable unitary connections on $E$ such that 
 $\Vert \Lambda F_{D_j}\Vert_{L^\infty}$ is bounded uniformly for all $j$.  Fix $p>4$.
Then there is:
\begin{enumerate}
\item   a subsequence $\{j_k\}$,
\item a  finite subset $\ansing\subset X$,
\item  a smooth hermitian vector bundle $(E_\infty, H_\infty)\to X\setminus \ansing$ with a finite action connection $D_\infty$ on $E_\infty$,
\item for any compact set $W\subset\subset X\setminus  \ansing$, an $L^p_2$-isometry 
$
\tau^W:(E_\infty, H_\infty)\bigr|_W \lra (E, H)\bigr|_W
$.
\end{enumerate}
such that for $W\subset W'\subset\subset X\setminus  \ansing$, $\tau^W=\tau^{W'}\bigr|_W$, and 
$
\tau^W(D_{j_k}) \weakarrow D_\infty
$ weakly\footnote{\, We will denote weak convergence by
``\, $\weakarrow$\, " and strong convergence by ``\, $\to$\, ".} in $L^p_1(W)$.
\end{Prop}

\noindent
We will call any  $D_\infty$ arising in this way an  \emph{Uhlenbeck limit} of the sequence $D_j$. 
 We will often omit the isometries $\tau^W$ from the notation, and simply identify $(E_\infty, H_\infty)$ with $(E,H)$ on $X\setminus \ansing$.  Also, it is useful to note here that weak $L^p_{1,loc.}$ convergence implies convergence of local holomorphic frames.  This may be proven, for example,  using Webster's proof of the Newlander-Nirenberg theorem \cite{W}.  We refer to \cite{DW} for more details on the following result:
 \begin{Prop}  \label{P:webster}
Let $D_j\weakarrow D_\infty$  in $L^p_{1, loc.}(X\setminus \ansing)$ for some $p>4$, as in Prop.\ \ref{P:compactness}. Then for each $x\in X\setminus \ansing$ there is:
\begin{enumerate}
\item  a coordinate neighborhood $U\subset X\setminus \ansing$ of $x$,
\item  a sequence $\{{\bf s}_j\}$ of $D_j^{\prime\prime}$-holomorphic frames on $U$,
\item a $D_\infty^{\prime\prime}$-holomorphic frame ${\bf s}_\infty$ on $U$,
\item  and a subsequence $\{ j_k\}\subset\{j\}$,
\end{enumerate}
 such that 
${\bf s}_{j_k}\to {\bf s}_\infty$ in $C^{1}(U)$.
\end{Prop}

 Next, we turn to a situation where the Uhlenbeck limits are Yang-Mills:

\begin{Prop} \label{P:ymcompactness}
If in addition to the assumptions in Prop.\ \ref{P:compactness} we assume
$
\Vert D_j\Lambda F_{D_j}\Vert_{L^2}\to 0
$,
then any Uhlenbeck limit $D_\infty$ is Yang-Mills.    Moreover, 
the triple $(E_\infty, D_\infty, H_\infty)$ extends smoothly to $X$, and
the extension has a holomorphic orthogonal splitting as a direct sum:
 $$\bigoplus_{i=1}^\ell (Q_{\infty}^{(i)}, D_{\infty}^{(i)}, H_{\infty}^{(i)})\ ,$$
where $H_{\infty}^{(i)}$ is a 
Hermitian-Einstein metric on $Q_{\infty}^{(i)}$.
\end{Prop}

\begin{proof}  The last statement follows by the removable singularities theorem (cf.\ \cite{U2}) and the argument cited above (cf.\ Prop.\ \ref{P:split}).  To see that $D_\infty$ is Yang-Mills, we argue as follows:  by the compactness of the embedding $L^p_1\hookrightarrow C^0$ and the fact that $D_{j_k}
\weakarrow  D_\infty$ weakly in $L^p_{1,loc.}$,  we may assume
$D_{j_k}  \overset{\scriptscriptstyle C^0_{loc.}}{\lra}  D_\infty$, and
$\Lambda F_{ D_{j_k}}\weakarrow\Lambda F_{ D_\infty}$ in $L^p_{loc.}$.
It follows that $D_\infty\Lambda F_{D_{j_k}}
 \weakarrow
 D_\infty\Lambda F_{D_\infty}$ in $L^2_{-1,loc.}$, say, where $L^2_{-k}$ is the dual space to $L^2_k$.
On the other hand, $
 D_\infty\Lambda F_{ D_{j_k}}=
 D_{j_k}\Lambda F_{ D_{j_k}}
+[D_\infty- D_{j_k},\Lambda F_{ D_{j_k}}]$,
so by the added hypothesis we also have,
$
 D_\infty\Lambda F_{ D_{j_k}}\overset{\scriptscriptstyle L^2_{loc.}}{\lra}  0$.
This  implies $ D_\infty\Lambda F_{ D_\infty}=0$.  The statement that $D_\infty$ is
Yang-Mills now follows from the K\"ahler identities (\ref{E:kahleridentities}).  This completes the proof.
\end{proof}

\begin{Cor} \label{C:lpconvergence}
With the assumptions as in Prop.\ \ref{P:ymcompactness}, $\Lambda F_{ D_{j_k}}\overset{\scriptscriptstyle L^p}{\lra} \Lambda F_{ D_\infty}$ for all $1\leq p<\infty$.
\end{Cor}

\begin{proof}
Set $f_k=\Lambda F_{ D_{j_k}}- \Lambda F_{ D_\infty}$.  Then by the proof of Prop.\ \ref{P:ymcompactness},
$f_k\weakarrow 0$ in $L^p_{loc.}$ and $D_\infty f_k\to 0$ strongly in $L^2_{loc.}$.  By Kato's inequality, $|f_k|$ is uniformly bounded in $L^2_{1,loc.}$, so $|f_k|\overset{\scriptscriptstyle L^2_{loc.}}{\lra} 0$.  Since $|f_k|$ is also uniformly bounded in $L^\infty$, it follows that $|f_k|\overset{\scriptscriptstyle L^p}{\lra} 0$  for all $p$.
\end{proof}

We conclude this subsection with a technical result on the boundedness of second fundamental forms.  This will be important in the proof of the main result. Let $\Omega\subset X$ be an open set.  Let  $\{D_j\}$, $D_\infty$ be integrable unitary connections on $E\to\Omega$ with $D_j\weakarrow D_\infty$ in $L^p_{1,loc.}(\Omega)$, for all $p$, $1\leq p<\infty$.  Let $\pi_j$ (resp.\ $\pi_\infty$) be (smooth) projections onto $D_j^{\prime\prime}$ (resp.\ $D_\infty^{\prime\prime}$) holomorphic subbundles of $E$. Let us assume the following:
\begin{enumerate}
\item $\Lambda F_{D_j}$ is  bounded in $L^\infty_{loc.}(\Omega)$ uniformly in $j$;
\item $\pi_j$ is bounded in $L^2_{1,loc.}(\Omega)$ uniformly in $j$;
\item $\Vert \pi_j-\pi_\infty\Vert_{L^\infty_{loc.}(\Omega)}\to 0$ as $j\to \infty$.
\end{enumerate}

\begin{Lem} \label{L:fundamentalform}
With the assumptions above, $\pi_j$ is bounded in $L^p_{2,loc.}(\Omega)$ uniformly in $j$, for all $p$.  In particular, the second fundamental forms $D_j^{\prime\prime}\pi_j$ are locally uniformly bounded.
 \end{Lem}
 
 \begin{proof}    
 Set $p_j=\pi_j-\pi_\infty$.  By \cite[Lemma 3.2]{D} we may write:
  \begin{equation} \label{E:bootstrap}
   \Delta_{D_\infty}(p_j)
  =\left\{D_\infty p_j,   D_\infty p_j\right\}+ \left\{ D_\infty \Gamma_j, p_j\right\} +  \left\{ D_\infty p_j , \Gamma_j\right\} + \left\{\Gamma_j,\Gamma_j\right\}+ G_j\ ,
  \end{equation} 
where   the brackets $\{\, , \}$  indicate a bilinear combination of the two arguments with bounded coefficients,
  $\Gamma_j=D_j-D_\infty$, and $G_j$ is uniformly bounded in $L^\infty_{loc.}(\Omega)$. Since $\Gamma_j$ is 
 uniformly bounded in $L^p_{1, loc.}(\Omega)$ and $p_j\to 0$ in $L^\infty_{loc.}(\Omega)$, it follows from \cite[Thm.\ 1.4]{GiM} and \cite[Thm. VI.1.5]{Gi} that $p_j$ is bounded in $C^{1,\alpha}_{loc.}(\Omega)$ for any $0<\alpha<1$,  uniformly in $j$.  The $L^p_{2,loc.}(\Omega)$ bound then follows from the  $L^p$-elliptic estimate (cf.\ \cite[Thm.\ 9.11]{GT}) applied to \eqref{E:bootstrap}. 
 For the second statement, write: $D_j^{\prime\prime}\pi_j=D_\infty^{\prime\prime}\pi_j +\Gamma_j^{\prime\prime}\pi_j$, and note that the right hand side is locally bounded since it is in $L^p_{1, loc.}$ for $p>4$.
 \end{proof}
 
\subsection{The Yang-Mills Flow}  \label{S:ymflow}

The basic object of interest in this paper is the   Yang-Mills
flow for a family of unitary connections $D=D(t)=D_t$.  This is the $L^2$-gradient flow of the YM functional, which  may be
written as follows: 
\begin{equation} \label{E:ymflow}
\frac{\partial D}{\partial t }= -D^\ast F_D \ ,\qquad
D(0)=D_0\in\Aone_H\ .
\end{equation}

\noindent
By the work of Donaldson and Simpson (cf.\ \cite{Do1, Si}), (\ref{E:ymflow}) has a unique solution in $(\Aone_H/\G)\times[0,\infty)$.  Furthermore, $D_t$ lies
in a single
$\GC$ orbit, namely $\GC\cdot D_0$,  for all $t\in[0,\infty)$.  One way to see this is to fix the $\dbar$-operator $\dbar_E=D_0^{\prime\prime}$ on $E$
and look at the family of hermitian metrics $H=H(t)=H_t$ satisfying the  Hermitian-Yang-Mills flow equations:
\begin{equation} \label{E:hymflow}
H_t^{-1}\frac{\partial H_t}{\partial t} = -2\left( \sqrt{-1}\, \Lambda_\omega F_{H_t}-\mu{\bf I}_E  \right) \ ,\qquad
H(0)=H_0\ .
\end{equation}
In the equation above, $F_{H_t}$ denotes the curvature of $D_t=(\dbar_E,H_t)$,
and $\mu=\mu_\omega(E)$ (see (\ref{E:slope})) depends only on the topology of $E$ and the K\"ahler form $\omega$. 
The two systems \eqref{E:ymflow} and \eqref{E:hymflow} are equivalent up to gauge.  See
\cite{Do1} for more details. Also notice that it is easy to factor out the trace part of the connection and gauge
transformations.  Therefore, in the following we shall assume that \emph{the solutions
to the above equations all preserve determinants}.  Also, in the following and throughout the paper, 
  we will often omit the K\"ahler form  $\omega$ from the notation when doing so is unambiguous.

The following  is an immediate consequence of (\ref{E:ymflow}) and (\ref{E:hymflow}) (see \cite[Prop.\ 16]{Do1} for a proof):
\begin{Lem} \label{L:subsolution}
\begin{enumerate}
\item  Let $D_t$ be a solution to (\ref{E:ymflow}).  Then 
$
\partial F_{D_t}/\partial t=-\Delta_t F_{D_t}$, and:
$$\frac{d}{dt}\Vert F_{D_t}\Vert^2_{L^2}=-2\Vert D_t^\ast F_{D_t}\Vert^2_{L^2}\leq 0
\ . $$
Hence, $t\mapsto \YM(D_t)$ and $t\mapsto \HYM(D_t)$ are nonincreasing.
\item  The pointwise norm $|\Lambda F_{D_t}|^2$ satisifies
$\displaystyle
\frac{\partial}{\partial t}|\Lambda F_{D_t}|^2+\Delta |\Lambda F_{D_t}|^2\leq 0
$.
\end{enumerate}
\end{Lem}
 
Uhlenbeck compactness applied to the flow gives the following:

\begin{Prop} \label{P:convergence}
 Let $D_t$ be a solution to \eqref{E:ymflow} on a compact K\"ahler surface $X$,
 and fix $p>4$.
For any sequence $t_j\to \infty$ we can find the following:
\begin{enumerate}
\item a subsequence $\{t_{j_k}\}$, 
\item a finite set of points $\ansing\subset X$,
\item  a smooth hermitian vector bundle $(E_\infty, H_\infty)\to X$ with a finite action Yang-Mills connection $D_\infty$ on $E_\infty$,
\item and on any compact set $W\subset\subset X\setminus \ansing$, an $L^p_2$-isometry 
$
\tau^W:(E_\infty, H_\infty)\bigr|_W \lra (E_0, H_0)\bigr|_W
$,
 \end{enumerate}
such that for $W\subset W'\subset\subset X\setminus \ansing$, $\tau^W=g^{W'}\bigr|_W$, and 
$
\tau^W(D_{t_{j_k}}) \weakarrow D_\infty,
$ weakly in $L^p_1(W)$.
Moreover, the triple $(E_\infty, D_\infty, H_\infty)$ extends smoothly to $X$, and the extension has a holomorphic orthogonal splitting as a
direct sum: $\oplus_{i=1}^\ell (Q_{\infty}^{(i)}, D_{\infty}^{(i)}, H_{\infty}^{(i)})$,
 where $H_{\infty}^{(i)}$ is a Hermitian-Einstein metric on $Q_{\infty}^{(i)}$.
\end{Prop}

\begin{proof}  
 By Lemma \ref{L:subsolution} (2) and  the maximum
principle,
$\Vert \Lambda F_{D_{t}}\Vert_{L^\infty}$ is decreasing in $t$, and is therefore uniformly bounded.
 By \cite[Prop.\ 6.2.14]{DoKr}, $\lim_{t\to\infty}\Vert D_t\Lambda F_{D_t}\Vert_{L^2}=0$.
%By Lemma \ref{L:subsolution} (1), we have:
%$$
%\int_0^\infty \Vert D_t^\ast F_{D_t}\Vert^2_{L^2}\, dt < +\infty\ .
%$$
%Hence, there is a sequence $t_j\to +\infty$ along which $\Vert D_{t_j}^\ast F_{D_{t_j}}\Vert_{L^2}\to 0$.  
 The weak $L^p_{1,loc.}$ convergence along a subsequence
to a YM connection now follows from Prop.\ \ref{P:ymcompactness}.
\end{proof}

As in \S \ref{S:ym}, we will call $D_\infty$ an Uhlenbeck limit of the flow.   Note that a priori,  $D_\infty$ may
depend on the choice of subsequence $\{t_{j_k}\}$, however we will see shortly that this is not the case.

\begin{Lem} \label{L:limit}
Let $D_{t_j}$ be a sequence of connections along the YM flow with Uhlenbeck limit $D_\infty$.  Then for $t_j\geq t_0 \geq 0$,
$
\Vert \Lambda F_{D_\infty}\Vert_{L^\infty}\leq \Vert  \Lambda F_{D_{t_j}}\Vert_{L^\infty}\leq \Vert
 \Lambda F_{D_{t_0}}\Vert_{L^\infty}
$.
\end{Lem}

\begin{proof}  As stated above,
$\Vert \Lambda F_{D_{t}}\Vert_{L^\infty}$ is decreasing in $t$.  Fix $t\geq 0$.   Then for any $1\leq p<\infty$ and $j$ sufficiently large we have:
 $
\Vert  \Lambda F_{D_{t_j}}\Vert_{L^p}  
\leq (2\pi)^{1/p}\Vert  \Lambda F_{D_{t_j}}\Vert_{L^\infty} 
\leq (2\pi)^{1/p} \Vert   \Lambda F_{D_{t}}\Vert_{L^\infty}   
$
(recall $\vol(X)=2\pi$).  On the other hand, by Cor.\ \ref{C:lpconvergence},
$
\lim_{j\to \infty}\Vert   \Lambda F_{D_{t_j}}\Vert_{L^p}  
=\Vert   \Lambda F_{D_{\infty}}\Vert_{L^p}  
 $,
 for all $p$.
 Hence, 
 $
 \Vert   \Lambda F_{D_{\infty}}\Vert_{L^p}
 \leq (2\pi)^{1/p} \Vert  \Lambda F_{D_{t}}\Vert_{L^\infty}  
 $.
 Letting $p\to \infty$, we conclude
 $
 \Vert \Lambda F_{D_{\infty}}\Vert_{L^\infty}
 \leq \Vert  \Lambda F_{D_{t}}\Vert_{L^\infty} 
 $.
   \end{proof}

\begin{Lem} \label{L:lp}
If $D_\infty$ is the Uhlenbeck limit of $D_{t_j}$, then $\Lambda F_{D_{t_j}}\overset{\scriptscriptstyle L^p}{\lra} \Lambda F_{D_\infty}$ for all $1\leq p<\infty$. 
Moreover, $\lim_{t\to \infty} \HYM(D_t)=\HYM(D_\infty)$.
\end{Lem}

\begin{proof} The first part of the lemma follows from Cor.\ \ref{C:lpconvergence}.   The second part is immediate, since by Lemma 
\ref{L:subsolution}, $t\mapsto\HYM(D_t)$ is nonincreasing, and
$\HYM(D_{t_j})\to\HYM(D_{\infty})$.
\end{proof}
  
 The Uhlenbeck limits obtained from the Yang-Mills flow are unique:
 
 \begin{Prop} \label{P:uniqueness}
 Let $D_t$  be the solution to the YM flow (\ref{E:ymflow}), and suppose $D_\infty$ is an Uhlenbeck limit for some sequence $D_{t_j}$ with singular set $\ansing$.  Then $D_t\to D_\infty$ in $L^2_{loc.}$ away from $\ansing$.  In particular, the Uhlenbeck limit of the flow is uniquely defined up to gauge.
 \end{Prop}  
 
 \begin{proof}
 By (\ref{E:chernclass}) and Lemma \ref{L:lp}, 
 $\displaystyle
 \lim_{t\to\infty} \YM(D_{t})=\HYM(D_\infty)+ 4\pi^2(2c_2(E)-c_1^2(E))
 $.
  From (\ref{E:ymflow}) and Lemma \ref{L:subsolution} (1), if $t_j\geq t$:
\begin{align*} 
\Vert D_{t_j}-D_t\Vert_{L^2}^2 &\leq \int_t^{t_j} \left\Vert \frac{\partial D_{s}}{\partial s}\right  \Vert_{L^2}^2\,       ds         
= \int_t^{t_j} \Vert   D_{s}^\ast  F_{D_{s}} \Vert_{L^2}^2\,    ds       \\
&= -\frac{1}{2}\int_t^{t_j}\frac{d}{ds} \Vert    F_{D_{s}} \Vert_{L^2}^2 \,    ds         
=\frac{1}{2}\left(\YM(D_t)-\YM(D_{t_j}) \right)\ . 
\end{align*}
  Since the limit of $\YM(D_t)$ as $t\to\infty$ exists, 
$\Vert D_{t_j}-D_t\Vert_{L^2}\to 0$ as $t$ and $j\to \infty$.  Since $D_{t_j}\overset{\scriptscriptstyle L^2_{loc.}}{\lra} D_\infty$, the convergence follows.    For the last statement, if $D_\infty$ and $\widetilde D_\infty$ are two Uhlenbeck limits with singular sets $\ansing$ and $ \ansingtilde$, then the argument above shows that $D_\infty$ and $\widetilde D_\infty$ are gauge equivalent on $X\setminus \ansing\cup \ansingtilde $.   In particular, the holomorphic bundles $(E_\infty, D_\infty^{\prime\prime})$ and $(\widetilde E_\infty,\widetilde D_\infty^{\prime\prime})$ are isomorphic on $X\setminus \ansing\cup\ansingtilde $.  But then their reflexive extensions are isomorphic as well by Hartogs theorem.
 \end{proof}

  \begin{Rem}  \label{R:limit}
  In light of Prop.\ \ref{P:uniqueness}, we may speak of \emph{the} Uhlenbeck limit of the flow $D_t$.  Note, however, that we have not established that  the singular set $\ansing$ is independent of the subsequence $\{t_j\}$.
  \end{Rem}

 We next turn to the HN type of the Uhlenbeck limit:
 
\begin{Lem} \label{L:degreesemicontinuity}
Let $D_j=g_j(D_0)$ be a sequence of complex gauge equivalent integrable connections on a complex vector bundle $E$ of rank $R$ with hermitian metric
$H_0$.  Let
$S$ be a coherent subsheaf of $(E, D_0^{\prime\prime})$ of rank $r$.  Suppose that 
  $\sqrt{-1}\,  \Lambda F_{D_j}\overset{\scriptscriptstyle L^1}{\lra} {\bf a}$, where ${\bf a}\in L^1(\sqrt{-1}\mathfrak{u}(E))$, and that the eigenvalues 
$\lambda_1\geq\cdots\geq \lambda_R$
of $\bf a$ (counted with multiplicities)  are constant almost everywhere.  Then:
$
\deg(S)\leq \sum_{i\leq r} \lambda_i
$.
\end{Lem}

\begin{proof}
Since $\deg(S)\leq\deg (\sat_E(S))$, we may assume that $S$ is saturated.
Let $\pi_j$ denote the orthogonal projection onto $g_j(S)$ with respect to the hermitian metric $H_0$. 
This is a bounded measurable hermitian endomorphism of $E$, smooth away from the singular set of $E/S$.  
 The condition of being a weakly holomorphic projection implies $\pi_j^2=\pi_j=\pi_j^\ast$, $\pi_j^\perp D_j^{\prime\prime}\pi_j=0$, where $\pi_j^\perp={\bf I}_E-\pi_j$. In fact, the $\pi_j$ are
$L^2_1$ sections of the smooth endomorphism bundle of $E$ (cf.\ \cite[\S 4]{UY}), and conversely, any such $\pi$ defines a unique saturated subsheaf\footnote{We will often confuse the notation $\pi$ of the projection operator with the subsheaf it defines.}.
 Moreover, the usual degree formula applies (see
\cite[Lemma 3.2]{Si}), so:
\begin{align}
\deg(S)&= \frac{1}{2\pi}\int_X\left(\tr \left(\sqrt{-1}\, \Lambda F_{D_j}\pi_j\right)-|D_j^{\prime\prime}\pi_j|^2\right)\, dvol
\leq \frac{1}{2\pi}\int_X \tr \left(\sqrt{-1}\, \Lambda F_{D_j}\pi_j\right)\, dvol \notag \\
&= \frac{1}{2\pi}\int_X \tr \left({\bf a}\pi_j\right)\ dvol +\frac{1}{2\pi}\int_X \tr \left(\left(\sqrt{-1}\, \Lambda F_{D_j}-{\bf
a}\right)\pi_j\right)\, dvol
\label{E:degbound}\ ,
\end{align}
 We now use the following result from linear algebra:
 \begin{Claim}
 Let $V$ be a finite dimensional hermitian vector space of complex dimension $R$, $L\in\End(V)$ a hermitian operator with eigenvalues
$\lambda_1\geq\cdots\geq\lambda_R$ (counted with multiplicities).  Let $\pi=\pi^2=\pi^\ast$ denote the orthogonal projection onto a
subspace of dimension $r$.  Then $\tr(L\pi)\leq\sum_{i\leq r}\lambda_i$.
\end{Claim} 
 \begin{proof}[Proof (sketch)]
 Let $\{e_i\}_{i=1}^R$ be a unitary basis with $Le_i=\lambda e_i$.  
 If we set   $\alpha_i=\Vert \pi e_i\Vert^2$,
then the claim follows by showing that the affine function
 $
 F(\alpha_1,\ldots,\alpha_R)=\sum_{i\leq r} \lambda_i-\sum_{i=1}^R \lambda_i\alpha_i 
 $, is nonnegative on the affine set $0\leq \alpha_i\leq 1$, $\sum_{i=1}^R \alpha_i=r$.
  This may be proven by considering the extreme values of the $\alpha_i$'s and using induction on $r$ and $R$.  We omit the details.
 \end{proof}
\noindent
 Given the claim, along with  the fact that $\Vert \pi_j\Vert_{L^\infty}\leq 1$, and the normalization $\vol(X)=2\pi$, 
it follows from (\ref{E:degbound}) that
  $
 \deg(S)\leq \sum_{i\leq r} \lambda_i +(1/2\pi)\left\Vert \sqrt{-1}\, \Lambda F_{D_j}-{\bf a}\right\Vert_{L^1} 
 $.
 Now let $j\to \infty$ in this inequality to complete the proof of the lemma.
\end{proof}

Recall the partial ordering (\ref{E:slopeordering}) of HN types of holomorphic structures on $E$:

\begin{Prop} \label{P:slopesemicontinuity}
Let $D_j$ be a sequence along the YM flow on a bundle $E$ of rank $R$  with Uhlenbeck limit $D_\infty$.  Let $\vec\mu_0=(\mu_1,\ldots, \mu_R)$ be the HN type of
$E$ with the holomorphic structure $D_0^{\prime\prime}$, and let
$\vec\lambda_\infty=(\lambda_1,\ldots,\lambda_R)$ the type of $D_{\infty}^{\prime\prime}$.  Then $\vec\mu_0\leq \vec\lambda_\infty$.
\end{Prop}

\begin{proof}
Let $\{E_i\}_{i=1}^\ell$  be the HN filtration  of $D_0^{\prime\prime}$.   Then
$
\deg(E_i)=\sum_{j\leq\rk(E_i)}\mu_j
$.
By Lemma \ref{L:lp}, $\Lambda
F_{D_j} 
 \overset {\scriptscriptstyle L^1} { \lra }\Lambda F_{D_\infty}$.
   The type
 $\vec\lambda_\infty$ corresponds to the (constant) eigenvalues of $\Lambda
F_{D_\infty}$.   Lemma \ref{L:degreesemicontinuity} applied to $S=E_i$ implies
$
\deg(E_i)\leq \sum_{j\leq\rk(E_i)}\lambda_j
$, for each $i=1,\ldots,\ell$.  The proposition now follows from Lemma \ref{L:shatz}.
  \end{proof}
 
 The following generalizes a result in \cite{AB} to K\"ahler surfaces:
 
 \begin{Cor}  \label{C:minimum}
  Let $\vec\mu=(\mu_1,\ldots, \mu_R)$ be the Harder-Narasimhan type of a rank $R$ holomorphic vector bundle $(E, \dbar_E)$ on  $X$.
  Then  for all  unitary connections $D$ in the $\GC$ orbit of $(E, \dbar_E)$:  
  $$
   \sum_{i=1}^R\mu_i^2 \leq\frac{1}{2\pi} \int_X |\Lambda F_{D}|^2 \, dvol\  . 
 $$
   \end{Cor}
 
 \begin{proof}
Let $D_t$ denote the YM flow with initial condition $D$. 
By Lemma \ref{L:subsolution} (2):
\begin{equation}  \label{E:decreasing}
 \int_X |\Lambda F_{D_t}|^2\, dvol \leq  \int_X |\Lambda F_{D}|^2\, dvol \ ,
\end{equation}
for all $t\geq 0$.  Let $D_\infty$ be the Uhlenbeck limit along a subsequence $t_j\to \infty$.  By Lemma
\ref{L:lp}:
\begin{equation}  \label{E:convergence}
 \int_X |\Lambda F_{D_\infty}|^2\, dvol =\lim_{j\to\infty} \int_X |\Lambda F_{D_{t_j}}|^2\, dvol \ .
\end{equation}
Now $\sqrt{-1}\Lambda F_{D_\infty}$ has constant eigenvalues $\vec\lambda_{\infty}$ which satisfy $\vec\mu\leq \vec\lambda_\infty$ by Prop.\ \ref{P:slopesemicontinuity}.  It follows from \cite[12.6]{AB} (see also Prop.\ \ref{P:key1} below) that
$
\sum_{i=1}^R\mu_i^2\leq  \sum_{i=1}^R\lambda_i^2
$.
This fact, along with (\ref{E:decreasing}) and (\ref{E:convergence}), prove the result.
\end{proof}

 %%%%%%%%%%%%%%%%%%%%%%%%%%%%%%%%%%%%%%%%%%%%%%%%%%%%%%%%%%%%%%%
  \subsection{Other Hermitian-Yang-Mills Type Functionals}  \label{S:other}
 
 One of the technical difficulties in dealing with holomorphic vector bundles of rank bigger than two is that the Hermitian-Yang-Mills
numbers do not distinguish the different critical levels (or equivalently, different Harder-Narasimhan types) of the functionals $\YM$ and $\HYM$.  This was resolved by Atiyah and Bott in the case of vector bundles over Riemann surfaces by introducing Yang-Mills type functionals corresponding to higher
symmetric functions of the eigenvalues of $\sqrt{-1}\, \ast F_D$ (cf.\   \cite[\S 8]{AB}).  In the case of vector bundles over higher dimensional K\"ahler
manifolds there are analytic restrictions on the type of functionals we may consider.  In this subsection we will explain in some
detail how to use these functionals in order to distinguish the different critical levels. 
 
 Let ${\mathfrak u}(R)$ denote the Lie algebra of the unitary group
$U(R)$. Fix a real number $\alpha \geq 1$.  Then for
${\bf a}\in{\mathfrak u}(R)$, a skew hermitian matrix with eigenvalues $\sqrt{-1}\, \lambda_1,\ldots,\sqrt{-1}\, \lambda_R$, let
$
\varphi_\alpha({\bf a})=\sum_{j=1}^R |\lambda_j|^\alpha
$.
 It is easy to see that we can find a family 
$\varphi_{\alpha, \rho}$, $0<\rho\leq 1$,  of smooth convex ad-invariant functions such that $\varphi_{\alpha,\rho}\to\varphi_\alpha$  uniformly on compact
subsets of
${\mathfrak u}(R)$ as $\rho\to 0$. Hence, 
by  \cite[Prop.\ 12.16]{AB} it follows that  $\varphi_\alpha$ is a convex function on ${\mathfrak
u}(R)$.  For a given number $N$, define:
\begin{equation} \label{E:hymalpha}
\HYM_{\alpha,N}(D)=\int_X\varphi_\alpha(\Lambda F_D+\sqrt{-1}N\, {\bf I}_E)\, dvol\ ,
\end{equation}
 and $\HYM_\alpha(D)=\HYM_{\alpha,0}(D)$.
Notice that $\HYM=\HYM_2$ is the ordinary HYM functional.  Also, by a slight abuse of notation, we will set
 $\HYM_{\alpha,N}\left(\vec\mu\right)=\HYM_\alpha(\vec\mu+N)=2\pi\varphi_\alpha(\sqrt{-1}\, (\vec\mu+N))$, where $\vec\mu+N=(\mu_1+N,\ldots,\mu_R+N)$ is identified
with  the diagonal matrix  $\diag(\mu_1+N,\ldots,\mu_R+N)$.
 In particular: 
\begin{equation} \label{E:hymminimum}
 \HYM\left( \vec\mu\right)=2\pi\sum_{i=1}^R \mu_i^2\ .
 \end{equation} 
\begin{Lem}  \label{L:norm}
The functional
$\displaystyle
{\bf a}\mapsto\left( \int_X \varphi_\alpha({\bf a})\ dvol\right)^{1/\alpha}
$,
defines a norm on $L^\alpha({\mathfrak u}(E))$ which is equivalent to the $L^\alpha$ norm. 
\end{Lem}

\begin{proof}
First, notice that there are universal constants $C,C'$ (depending on $R$) such that for any real numbers $\lambda_1,\ldots,\lambda_R$, and $\alpha\geq 1$:
$$
\frac{1}{C}\left(\sum_{i=1}^R|\lambda_i|^2\right)^{\alpha/2}
\leq 
\frac{1}{C}\left(\sum_{i=1}^R|\lambda_i|\right)^{\alpha}
\leq 
\sum_{i=1}^R|\lambda_i|^\alpha
\leq
C\left(\sum_{i=1}^R|\lambda_i|\right)^{\alpha}
\leq
C'\left(\sum_{i=1}^R|\lambda_i|^2\right)^{\alpha/2}\ .
$$
Applying this to the eigenvalues of ${\bf a}$ and by integrating over $X$:
$$
\frac{1}{C}\int_X\left(\tr {\bf a}{\bf a}^\ast\right)^{\alpha/2} dvol\leq \int_X\varphi_\alpha({\bf a})\, dvol\leq C'\int_X\left(\tr
{\bf a}{\bf a}^\ast\right)^{\alpha/2}\, dvol\ .
$$
The lemma follows.
\end{proof}

We will require three  important properties of the functionals $\HYM_{\alpha,N}$:
 
\begin{Prop} \label{P:key1}
 \begin{enumerate}
 \item
If $ \vec\mu\leq \vec\lambda$, then $\varphi_\alpha(\sqrt{-1}\, \vec\mu)\leq\varphi_\alpha(\sqrt{-1}\, \vec\lambda)$ for all $\alpha\geq 1$.  
 \item Assume $\mu_R\geq 0$ and $\lambda_R\geq 0$.  If
$\varphi_\alpha(\sqrt{-1}\, \vec\mu)=\varphi_\alpha(\sqrt{-1}\, \vec\lambda)$ for all $\alpha$ in some set
$A\subset[1,\infty)$ possessing a limit point, then $ \vec\mu= \vec\lambda$.
 \end{enumerate}
\end{Prop}

\begin{proof}
(1) follows from \cite[12.6]{AB}.  For (2), consider $f(\alpha)=\varphi_\alpha(\sqrt{-1}\, \vec\lambda)$ and
$g(\alpha)=\varphi_\alpha(\sqrt{-1}\, \vec\mu)$ as  functions of
$\alpha$.   As complex valued functions,  $f,g$ clearly have  analytic extensions to $\CBbb\setminus\{\alpha\leq 0\}$.  Suppose  that
$f(\alpha)=g(\alpha)$ for all
$\alpha \in A$.  Then by analyticity, $f(\alpha)=g(\alpha)$ for all $\alpha\in\CBbb\setminus\{\alpha\leq 0\}$.  If $\vec\mu\neq\vec\lambda$, then there is some $k$, $1\leq k\leq R$, such that $\mu_i=\lambda_i$ for $ i<k$, and  $\mu_k\neq \lambda_k$; say, $\mu_k>\lambda_k$ .   Then for any $\alpha>0$:
$$
\left(\frac{\mu_k}{\lambda_k}\right)^\alpha\leq \sum_{i=k}^R\left(\frac{\mu_i}{\lambda_k}\right)^\alpha
=\sum_{i=k}^R\left(\frac{\lambda_i}{\lambda_k}\right)^\alpha\leq R\ .
$$
Letting $\alpha\to\infty$, we obtain a contradiction; hence, the result.
\end{proof} 
 
 \begin{Prop} \label{P:key2}
Let $D_t$ be a solution of (\ref{E:ymflow}). Then for any $\alpha \geq 1$ and any $N$, $t\mapsto \HYM_{\alpha,N}(D_t)$ is nonincreasing.
\end{Prop}

\begin{proof}  Because we can approximate $\varphi_\alpha$ by  smooth convex ad-invariant functions $\varphi_{\alpha,\rho}\to\varphi_\alpha$, it suffices to show that the functional
$
t\mapsto\int_X\varphi_{\alpha,\rho}(\Lambda F_{D_t}+\sqrt{-1}N\, {\bf I}_E)\, dvol
$,
is nonincreasing along the flow for any $\rho>0$.  This follows from integrating the following inequality:
\begin{equation} \label{E:subphi}
(\partial/\partial t)\varphi_{\alpha,\rho}\left(\Lambda F_{D_t}+\sqrt{-1}N\, {\bf I}_E\right)+\Delta\varphi_{\alpha,\rho}\left(\Lambda F_{D_t}+\sqrt{-1}N\, {\bf I}_E\right)\leq 0\ .
\end{equation}
To prove (\ref{E:subphi}), simplify the notation by  setting $f =\Lambda F_{D_t}+\sqrt{-1}N\, {\bf I}_E$,
$\varphi=\varphi_{\alpha,\rho}$.  We first claim that:
\begin{equation} \label{E:laplacian}
\Delta(\varphi\circ f)(x) =-\ast\varphi^{\prime\prime}_{f(x)}\left( \ast D_t f, D_t f\right)+\varphi^{\prime}_{ f(x)}\left(\Delta_{D_t} f\right)\ .
\end{equation}
Indeed,
$
\Delta(\varphi\circ f)(x) =-\ast d\ast d(\varphi\circ f)(x)=-\ast d\ast\varphi^{\prime}_{ f(x)}(D_t  f)
$,
because $\varphi$ is invariant under the adjoint action.  Then:
\begin{align*}
-\ast d\ast\varphi^{\prime}_{ f(x)}(D_t  f)&=
-\ast d\left(\varphi^{\prime}_{ f(x)}(\ast D_t  f)\right) 
= -\ast\varphi^{\prime\prime}_{ f(x)}\left(\ast D_t f, D_t f\right)-\ast \varphi^{\prime}_{ f(x)}(D_t\ast D_t  f) \\
&= -\ast\varphi^{\prime\prime}_{ f(x)}\left(\ast D_t f, D_t f\right)+\varphi^{\prime}_{ f(x)}\left(\Delta_{D_t}  f\right) \ .
\end{align*}
The claim (\ref{E:laplacian}) follows.    Since $\varphi^{\prime\prime}$ is a positive definite quadratic form:
\begin{align*}
\Delta\left(\varphi\circ(\Lambda F_{H_t}+\sqrt{-1}N\, {\bf I}_E)\right)(x)&\leq \varphi^{\prime}\left(\Delta_{D_t} (\Lambda F_{D_t}+\sqrt{-1}N\, {\bf I}_E)\right)\\
&=  -\varphi^{\prime}\left(\frac{\partial}{\partial t}(\Lambda F_{D_t}+\sqrt{-1}N\, {\bf I}_E)\right)\\
 &=
-\frac{\partial}{\partial t}\varphi(\Lambda F_{D_t}+\sqrt{-1}N\, {\bf I}_E)\ ,
\end{align*}
where the first equality follows from Lemma \ref{L:subsolution} (1).  This verifies (\ref{E:subphi}) and completes the proof.
 \end{proof}

\begin{Prop} \label{P:equality}
Let $D_\infty$ be a subsequential Uhlenbeck limit of $D_t$, where $D_t$ is a solution to (\ref{E:ymflow}).  Then for any $\alpha\geq 1$ and any $N$,
$
\lim_{t\to\infty} \HYM_{\alpha,N}(D_t)=\HYM_{\alpha,N}(D_\infty)
$.
\end{Prop}

\begin{proof}  Let $D_\infty$ be the Uhlenbeck limit of a sequence $D_{t_j}$.  By  Lemma \ref{L:lp}, $\Lambda F_{D_{t_j}}\overset{\scriptscriptstyle L^p}{\lra}\Lambda F_{D_\infty}$  for all $p$.  Hence, by Lemma \ref{L:norm}, $\HYM_{\alpha,N}(D_{t_j})\to\HYM_{\alpha,N}(D_\infty)$.  The convergence in general follows by Prop.\
\ref{P:key2}, since $\HYM_{\alpha,N}(D_t)$ is nonincreasing in $t$.
\end{proof}

%%%%%%%%%%%%%%%%%%%%%%%%%%%%%%%%%%%%%%%%%%%%%%%%%%%%%%%%%%%%%%%%%%%%%%%%%%%%%%%%%%%%%%%%%%%%%%%%%%%%%%%%%%%%%%%%%%
%%%%%%%%%%%%%%%%%%%%%%%%%%%%%%%%%%%%%%%%%%%%%%%%%%%%%%%%%%%%%%%%%%%%%%%%%%%%%%%%%%%%%%%%%%%%%%%%%%%%%%%%%%%%%%%%%%

\section{Blow-up of the Harder-Narasimhan Filtration and Approximate Metrics}       \label{S:blowup}

%%%%%%%%%%%%%%%%%%%%%%%%%%%%%%%%%%%%%%%%%%%%%%%%%%%%%%%%%%%%%%%%%%%%%%%%%%%%%%%%%%%%%%%%%%%%%%%%%%%%%%%%%%%%%%%%%
%%%%%%%%%%%%%%%%%%%%%%%%%%%%%%%%%%%%%%%%%%%%%%%%%%%%%%%%%%%%%%%%%%%%%%%%%%%%%%%%%%%%%%%%%%%%%%%%%%%%%%%%%%%%%%%%%%

The first goal of this section is to show how to resolve the Harder-Narasimhan filtration  of a holomorphic bundle $E\to X$
 by passing to a modification $\pi: \widehat X\to X$ of the
K\"ahler surface  $X$.  We shall see that this procedure works well when the associated graded object consists of stable sheaves.  
While this result is not directly needed for the remaining sections, we have chosen to present it here since it may be of independent
interest. In the case where semistable quotients appear, the relationship between the HN filtration of $E$ and that of
$\pi^\ast E$  is more complicated, and we have not attempted a complete description.  As pointed out by the referee, much of this analysis has already appeared in the work of Buchdahl \cite{Bu1,Bu2}.

The second goal is to show that there is      a Hermitian metric
$\widehat H$ on $\pi^\ast E$  so that the Hermitian-Yang-Mills numbers of $(\pi^\ast E, \widehat H)$ with respect to a natural family of K\"ahler
metrics
$\omega_\varepsilon$ are arbitrarily close (in an appropriate norm) to the slopes of
$\Grhn_\omega(E)$ on $X$ for all $\varepsilon$ sufficiently small.  This  is an important first step toward finding
an approximate critical hermitian structure  on $X$ itself.  We formulate this latter result below in Thm.\ \ref{T:approximate}.
The argument we give here circumvents the need for an explicit description of the
HN filtration of
$\pi^\ast E$.

%%%%%%%%%%%%%%%%%%%%%%%%%%%%%%%%%%%%%%%%%%%%%%%%%%%%%%%%%%%%%%%%%%%%%%%%%%%%%%%%%%%%%%%%%%%%%%%%%%%%%%%%%%%%%%%%%%

\subsection{Resolution of the Harder-Narasimhan Filtration}  \label{S:hnresolution}

 By Prop.\ \ref{P:hnfiltration}, a holomorphic bundle $E$ admits a filtration by saturated
subsheaves
$E_i$ so that the successive quotients $Q_i=E_i/E_{i+1}$ are semistable and torsion-free.  It follows  that the $E_i$ are locally free sheaves, i.e.\ vector bundles.  They may, however, fail to be
subbundles at finitely many points.  Equivalently, the quotients $Q_i$ are not necessarily locally free.  For each $i$, we have an exact
sequence of sheaves:
$
0\to Q_i\to Q_i^{\ast\ast}\to T_i\to 0
$,
where $Q_i^{\ast\ast}$ is locally free and $T_i$ is a torsion sheaf supported at finitely many points.
Define the set $Z_i$ to be the support of $T_i$, and let $\algsing=\cup_{i=1}^k Z_i$.  We will refer to $\algsing$ as the  singular set of
the filtration $\{E_i\}$ (in this paper we ignore multiplicities).  We will prove the following:

\begin{Thm} \label{T:hnresolution}
Let $E\to X$ be a holomorphic vector bundle over a smooth K\"ahler surface $X$ with K\"ahler form $\omega$.   Let  $\Grhn_\omega(E)=\bigoplus_{i=1}^k Q_i$ be the
$\omega$-HN filtration of
$E$.  We assume that the $Q_i$ are \emph{stable}.
Then there is a smooth surface $\widehat X$ obtained from $X$ by a sequence of monoidal
transformations  $\pi:\widehat X\to X$ with exceptional set ${\bf e}\subset\widehat X$ satisfying the following properties:
 \begin{enumerate}
\item  $\pi({\bf e})=\algsing$; 
\item   There exists a smooth, closed $(1,1)$ form $\eta$ on $\widehat X$ and a number $\varepsilon_0>0$ such that  $\omega_\varepsilon=\pi^\ast\omega+\varepsilon\eta$ is a family K\"ahler metrics on $\widehat X$ for all $\varepsilon_0\geq \varepsilon>0$.
 \item There is a number $\varepsilon_1$, $\varepsilon_0\geq \varepsilon_1>0$  such that for every $0<\varepsilon\leq \varepsilon_1$, the
$\omega_\varepsilon$-HN filtration $\left\{\filt_{i,\omega_\varepsilon}(\widehat E)\right\}$ of $\widehat E=\pi^\ast E$ is independent of
$\varepsilon$ and is a filtration by \emph{subbundles}.
\item For $0<\varepsilon\leq\varepsilon_1$, $\varepsilon_1$ as above, $\filt_{i,\omega_\varepsilon}(\widehat E)=\sat_{\widehat E}\left(\pi^\ast\filt_i(E)\right)$.
 Moreover, if we write:  $\Grhn_{\omega_\varepsilon}(\widehat E)=\bigoplus_{i=1}^{\hat k} \widehat
Q_i$, then for $0<\varepsilon\leq \varepsilon_1$, $\hat k=k$, and $(\pi_\ast\widehat Q_i)^{\ast\ast}\simeq Q_i^{\ast\ast}$.
\end{enumerate}
\end{Thm}
 
 \begin{Rem} \label{R:semistable}
 The assumption that $Q_i$ is stable is necessary.  One can find examples where the pull-back of a semistable bundle is unstable for all $\varepsilon>0$. More generally, the resolution of the Harder-Narasimhan filtration of $E$ may not coincide with the Harder-Narasimhan filtration of $\pi^\ast E$ for any $\varepsilon>0$.
 \end{Rem}

 We begin by comparing stability of sheaves $\widehat E\to\widehat X$ with
stability of their direct images $E=\pi_\ast \widehat E$ on $X$.   To do this, we need to define K\"ahler metrics.   First, consider the case where $\widehat X$ is the blow-up of $X$ at a point and $\bf e$ is the exceptional divisor.   Then there is a smooth,  closed  form $\eta$ of type $(1,1)$ in the class of $c_1\left({\mathcal O}_{\widehat X}(-{\bf e})\right)$ on $\widehat X$ such that $ \omega_\varepsilon =\pi^\ast\omega +\varepsilon\eta$ is positive for all $\varepsilon >0$ sufficiently small.  This can be constructed quite explicitly (cf.\ \cite[pp.\ 182-187]{GH}).  In general, since $\widehat X$ is a sequence of blow-ups at points, we can construct a family of K\"ahler forms on $\widehat X$ be iterating the above argument.  We state this precisely as:
  \begin{Lem} \label{L:kahler}
  Let $\pi:\widehat X\to X$ be a sequence of monoidal transformations with exceptional set ${\bf e}$, and choose a K\"ahler form $\omega$ on $X$.  Then there is a smooth, closed $(1,1)$ form $\eta$ on $\widehat X$ and a number $\varepsilon_0>0$  with the following properties:
  \begin{enumerate}
   \item  $\omega_\varepsilon =\pi^\ast\omega +\varepsilon\eta$  is a K\"ahler form on $\widehat X$  for all $\varepsilon_0\geq \varepsilon>0$;
     \item  For any closed 2-form $\alpha$ on $X$, $\displaystyle \int_{\widehat X} \pi^\ast\alpha\wedge\eta =0$.
 \end{enumerate} 
\end{Lem}
  
 Consider a  family of K\"ahler forms  in the manner of Lemma \ref{L:kahler}. Note that we do not normalize the volume of $(\widehat X,\omega_\varepsilon)$, though we still assume the normalization on $(X,\omega)$; so $\vol(\widehat X,\omega_\varepsilon)\to 2\pi$ as $\varepsilon\to 0$.
 In the following, let us agree that the slope $\mu(E)$ of a sheaf on $X$ will be taken with respect to $\omega$.  For a sheaf
$\widehat E$ on $\widehat X$, we denote by $\mu_{\varepsilon}(\widehat E)$ the slope of $\widehat E$ with respect
to the metric $\omega_\varepsilon$.  Similarly, a subscript $\varepsilon$ will indicate that the quantity in question is taken with respect to
$\omega_\varepsilon$.  With this understood,  we have the following:
\begin{Prop} \label{P:stability}
Given a holomorphic vector bundle $\widehat E\to\widehat X$ with $E=\pi_\ast \widehat E$, and given $\delta>0$,
there is $\varepsilon_1>0$, depending upon $\widehat E$, such that for all
$0<\varepsilon\leq
\varepsilon_1$ we have the following inequalities:
\begin{enumerate}
\item  $\mu(E)-\delta\leq \mu_\varepsilon(\widehat E)\leq \mu(E)+\delta$,
\item  $\mu_{max}(E)-\delta\leq \mu_{max,\varepsilon}(\widehat E)\leq \mu_{max}(E)+\delta$,
\item $\mu_{min}(E)-\delta\leq \mu_{min,\varepsilon}(\widehat E)\leq \mu_{min}(E)+\delta$.
\end{enumerate}
\end{Prop}

\begin{proof}
Since $\mu_{min}(E)=-\mu_{max}(E^\ast)$,
part (3) will follow from part (2) applied to $\widehat E^\ast$.  
Parts (1) and (2) are essentially contained in \cite[Lemma 5]{Bu1}.  The statement there assumed $\widehat E$ is a pull-back bundle, but the proof works as well for general $\widehat E$.
\end{proof}

As a consequence, we have (cf.\ \cite[Prop.\ 3.4 (d)]{Bu2}):

\begin{Cor} \label{C:stability}
Let $\widehat E\to\widehat X$ and $E=\pi_\ast \widehat E$ be as above.  If $E$ is $\omega$-stable, then there is
a number $\varepsilon_2>0$,  depending upon $\widehat E$, such that $\widehat
E$ is $\omega_\varepsilon$-stable for all $0<\varepsilon\leq \varepsilon_2$.
\end{Cor}

  An inductive argument repeatedly using Prop.\ \ref{P:stability} implies convergence of the HN type:
 
 \begin{Cor} \label{C:hntype}
Let  $\widehat E\to\widehat X$ be a holomorphic vector bundle with $E=\pi_\ast \widehat E$.  Let $\vec\mu_\varepsilon$ denote the HN type of $\widehat E$ with respect to $\omega_\varepsilon$ and 
 $\vec\mu$  the HN type of $ E$ with respect to $\omega$.  Then $\vec\mu_\varepsilon\to\vec\mu$ as $\varepsilon\to 0$.
 \end{Cor}

Next, we state a general result on resolution of filtrations:

\begin{Prop}  \label{P:resolution2}
Let
$0=E_0\subset E_1\subset\cdots\subset E_{\ell-1}\subset E_\ell=E$,
be a filtration of a holomorphic vector bundle $E\to X$ by saturated subsheaves $E_i$, and set $Q_i=E_i/E_{i-1}$.  Then there is a sequence of
monoidal transformations
$\pi : \widehat X\to X$ with exceptional set ${\bf e}$ and a filtration
$0=\widehat E_0\subset \widehat E_1\subset\cdots\subset \widehat E_{\ell-1}\subset \widehat E_\ell=\widehat E=\pi^\ast E$,
such that each $\widehat E_i=\sat_{\widehat E}(\pi^\ast E_i)$ is a \emph{ subbundle} of $\widehat E$.   If we let $\widehat Q_i=\widehat
E_i/\widehat E_{i-1}$, we also have  exact sequences
$0\to Q_i\to \pi_\ast \widehat Q_i\to T_i\to 0$,
where $T_i$ is a torsion sheaf supported at the singular set of $Q_i$.  Moreover, $\pi({\bf e})=\algsing$, the union of the 
singular sets of
$Q_i$;
$\pi_\ast \widehat E_i=E_i$; and $Q_i^{\ast\ast}=(\pi_\ast\widehat Q_i)^{\ast\ast}$.
\end{Prop}

\begin{proof}
 The proof is standard resolution of singularities (cf.\  \cite[\S 3]{Bu1} for the step 2 filtration; the general argument then follows by induction).  The form $\widehat E_i=\sat_{\widehat E}(\pi^\ast E_i)$ follows from  Lemma \ref{L:sat}.  The remaining statements follow easily, and we omit the details.
\end{proof}

\begin{Prop} \label{P:hnresolution}
  Let $\pi : \widehat X\to X$ be  a sequence of monoidal transformations   with exceptional set ${\bf e}$ as above.  Let
$\widehat E\to\widehat X$ be a holomorphic vector bundle, and set $E=\pi_\ast \widehat E$.  Let $E_i=\filt_i(E)$ denote the HN filtration of $E$,
and assume that the successive quotients $E_i/E_{i-1}$ are \emph{stable}.  Let $\hat\jmath:\pi^\ast E_i\to \widehat E$ denote the induced map.
 We also assume that the sheaves $\sat_{\widehat E}(\hat\jmath(\pi^\ast E_i))$ are
subbundles of $\widehat E$.
 Then for $\varepsilon>0$ sufficiently small, the HN filtration
$\{\filt_{i,\varepsilon}(\widehat E)\}$ with respect to the K\"ahler metrics of Lemma \ref{L:kahler} is independent of $\varepsilon$ and is given
by
$\filt_{i,\varepsilon}(\widehat E)=\sat_{\widehat E}(\hat\jmath(\pi^\ast E_i))$.
\end{Prop}

\begin{proof}
We will proceed by induction on the length of the HN filtration of $E$.  If $E$ is stable, then the result follows from Cor.\ 
\ref{C:stability}.  
Assume now that $E$ is not stable, and \emph{define} $\widehat E_i=\sat_{\widehat E}(\hat\jmath(\pi^\ast E_i))$. Note that with this definition it follows as in the proof of Prop.\ \ref{P:stability} that
$\pi_\ast\widehat E_i=E_i$. We claim that
$\widehat E_1=\filt_{1,\varepsilon}(\widehat E)$ for $\varepsilon$ sufficiently small.  This follows because (1) by the hypothesis and Cor.\
\ref{C:stability}, $\widehat E_1$ is stable for $\varepsilon$ sufficiently small, and (2) we may arrange that 
$\mu_\varepsilon(\widehat E_1)>\mu_{max}(\widehat Q_1)$ for $\varepsilon$ sufficiently small, where $\widehat Q_1=\widehat E/\widehat E_1$.  The
claim then follows from Prop.\ \ref{P:general} (3).  Let $Q_1=E/E_1$.  By pushing forward, we have
$
0\to Q_1\to \widetilde Q_1\to T\to 0
$,
where $\widetilde Q_1=\pi_\ast \widehat Q_1$, and $T$ is a torsion-sheaf supported at points.  Hence, by Prop.\  \ref{P:general} (1), the
HN filtrations of $Q_1$ and $\widetilde Q_1$ are related by
$
\filt_i(Q_1)=\ker(\filt_i(\widetilde Q_1)\to T)
$.
For convenience, set $F_i=\filt_i(Q_1)$, $\widetilde F_i=\filt_i(\widetilde Q_1)$.
Notice that $\widetilde Q_1$ and $\widehat Q_1$ continue to satisfy the hypothesis of the proposition.  Hence, by induction, we may assume that for $\varepsilon$ sufficiently small 
the HN filtration of $\widehat Q_1$ is given by $\filt_{i,\varepsilon}(\widehat Q_1)=\{ \sat_{\widehat Q_1}(\hat\jmath(\pi^\ast\widetilde F_i))\}$ (where $\hat\jmath$ now is the induced map to $\widehat Q_1$).
Now $\pi^\ast F_i\hookrightarrow \pi^\ast\widetilde F_i$ with a torsion quotient.  By Lemma \ref{L:sat},
$
\sat_{\widehat Q_1}(\hat\jmath(\pi^\ast\widetilde F_i))=\sat_{\widehat Q_1}(\hat\jmath(\pi^\ast F_i))
$.
This, combined with Prop.\ \ref{P:general} (2) gives
$
\filt_{i,\varepsilon}(\widehat E)=\ker(\widehat E\to\widehat Q_1\bigr/\sat_{\widehat Q_1}(\hat\jmath(\pi^\ast F_{i-1})))
$.
 Moreover, $F_{i-1}=E_i/E_1$.
Clearly then, $\filt_{i,\varepsilon}(\widehat E)$ contains $\hat\jmath(\pi^\ast E_i)$ with a torsion quotient.  Thus, again by Lemma \ref{L:sat}, 
$\filt_{i,\varepsilon}(\widehat E)=\sat_{\widehat E}(\hat\jmath(\pi^\ast E_i))$.
\end{proof}

\begin{proof}[Proof of Thm.\ \ref{T:hnresolution}]
By Prop.\ \ref{P:resolution2}, the HN filtration $\{E_i=\filt_{i,\omega}(E)\}$ of $E$ admits a resolution to a filtration $\{\widehat E_i\}$ by
subbundles on
$\widehat X$.  By Prop.\
\ref{P:hnresolution}, the filtration $\{\widehat E_i\}$ is the HN filtration with respect to $\omega_\varepsilon$ for $\varepsilon $ sufficiently
small.  The remaining assertions follow easily.
\end{proof}

%%%%%%%%%%%%%%%%%%%%%%%%%%%%%%%%%%%%%%%%%%%%%%%%%%%%%%%%%%%%%%%%%%%%%%%%%%%%%%%%%%%%%%%%%%%%%%%%%%%%%%%%%%%%%%%%%%%%%%%%%%%%%%%%%%%%%%%%%%%%%%%%%%

\subsection{An Approximate Critical Hermitian Structure} \label{S:approximate}

For a fixed holomorphic structure on $E\to X$, a critical point of the functional $H\mapsto \HYM(\dbar_E,H)$ is called a critical hermitian structure \cite[p.\ 108]{Ko}.  By the K\"ahler identities \eqref{E:kahleridentities}, this occurs if and only if the connection $(\dbar_E,H)$ is Yang-Mills.  
The general form for a critical hermitian structure  is therefore (see Prop.\
\ref{P:split}):
\begin{equation} \label{E:ym}
\sqrt{-1}\, \Lambda_\omega F_{(\dbar_E, H)}=
 \mu_1{\bf I}_{Q_1} \oplus\cdots\oplus \mu_\ell{\bf I}_{Q_\ell} \ .
\end{equation}
In the above, the holomorphic structure $\dbar_E$ on $E$ splits $\displaystyle E=\oplus_{i=1}^\ell Q_i$, and the induced metric on each factor
$Q_i$ is Hermitian-Einstein with slope $\mu_i$.  Notice that if we
assume the slopes are ordered
$\mu_1>\cdots>\mu_\ell$, then the HN filtration of $E$ is given by:
$
\filt_i(E)=\oplus_{j\leq i} Q_j
$.

For a general holomorphic structure on $E$, the HN filtration will not be holomorphically split, so there can exist no smooth metric satisfying
(\ref{E:ym}).  What is more, the HN filtration may not be given by subbundles, so the right hand side of (\ref{E:ym}) is not even  everywhere
defined as a smooth endomorphism.  In this subsection, we define precisely what is meant by an \emph{approximate} solution to (\ref{E:ym}) (compare the following
discussion  with that in
\cite[IV. \S5]{Ko}).

Let $H$ be a smooth metric on $E$, and let $\mathcal{F}=\{ F_i\}_{i=1}^\ell$ be a filtration of $E$ by saturated subsheaves:
$\mathcal{F} : 0=F_0\subset F_1\subset\cdots\subset F_\ell=E$.  Associated to each $F_i$ and the metric $H$ we have the unitary projection
$\pi_i^H$ onto $F_i$.  As mentioned previously, the $\pi_i^H$ are bounded $L^2_1$ hermitian endomorphisms.  For
convenience, we set
$\pi_0^H=0$.  
Next, suppose we are given a collection of real numbers $\mu_1,\ldots, \mu_\ell$.  From the data $\mathcal{F}$ and $(\mu_1,\ldots,\mu_\ell)$ we
define a bounded $L^2_1$ hermitian endomorphism of $E$ by 
$\Psi(\mathcal{F}, (\mu_1,\ldots,\mu_\ell),H)=\sum_{i=1}^\ell \mu_i(\pi_i^H-\pi_{i-1}^H)
$.
At points where the $F_i$ are all subbundles there is a smooth orthogonal splitting: 
$E=\oplus_{i=1}^\ell F_i/F_{i-1}$
 with respect to which
$\Psi(\mathcal{F}, (\mu_1,\ldots,\mu_\ell),H)$ is diagonal with entries $\mu_i$.
 Given a holomorphic hermitian vector bundle $E$ on a compact K\"ahler surface $(X,\omega)$, the  Harder-Narasimhan projection,  $\Psi_{\omega}^{hn}(\dbar_E,H)$, is the bounded $L^2_1$ hermitian endomorphism defined above in the particular case where $\mathcal F$ is the HN filtration $F_i=\filt_i(E)$ and $\mu_i=\mu(F_i/F_{i-1})$.

 \begin{Def} \label{D:approximate}
Fix $\delta >0$ and $1\leq p\leq \infty$.
An  $L^p$-$\delta$-approximate critical hermitian structure on  a holomorphic bundle $E$ is a smooth metric $H$ such that
$
\Vert\sqrt{-1}\, \Lambda_\omega F_{(\dbar_E, H)}-\Psi_{\omega}^{hn}(\dbar_E, H)\Vert_{L^p(\omega)}\leq \delta
$.
\end{Def}
 
Let us immediately point out the following:

\begin{Thm} \label{T:bundlecase}
If the HN filtration of $E$ is given by subbundles, then for any $\delta >0$ there is an $L^\infty$-$\delta$-approximate critical hermitian structure on
$E$.
\end{Thm}

\begin{proof}
First,   by the equivalence of holomorphic structures and integrable unitary connections, it suffices to show that for a \emph{fixed}
hermitian metric $H$ there is a smooth complex gauge transformation $g$ preserving the HN filtration such that:
\begin{equation} \label{E:gauge}
\left\Vert\sqrt{-1}\, \Lambda_\omega F_{(g(\dbar_E), H)}-\Psi_{\omega}^{hn}(g(\dbar_E), H)\right\Vert_{L^\infty(\omega)}\leq \delta\ ,
\end{equation}
(see \cite{Do1}).
Next, for semistable $E$ (i.e.\ the length $1$ case), the result follows by the
convergence 
$
\left\Vert\sqrt{-1}\, \Lambda_\omega F_{D_t}-\mu(E){\bf I}_E\right\Vert_{L^\infty(\omega)}\to 0
$,
where $D_t$ is a solution to the Yang-Mills flow equations (\ref{E:ymflow}) with any initial condition (cf.\ \cite[Cor.\ 25]{Do1}).   With this understood,  choose $\delta'$-approximate metrics, where $0<\delta'<<\delta$, on the semistable quotients $Q_i$ of the HN filtration of $E$ to fix a metric $H$ on $E=Q_1\oplus \cdots\oplus Q_\ell$.  Then by appropriately scaling the extension classes: $0\to E_{i-1}\to E_i\to Q_i\to 0$, one finds a complex gauge transformation satisfying (\ref{E:gauge}). We omit the details.
\end{proof}

We may now formulate one of the main results of this paper:

\begin{Thm} \label{T:approximate}
Let $E\to (X,\omega)$ be a holomorphic vector bundle on a smooth K\"ahler surface $X$.  Given any $\delta >0$ and any  $1\leq p<\infty$,
there is an $L^p$-$\delta$-approximate critical hermitian structure on $E$.
\end{Thm}

\begin{Rem} 
The metric produced in Thm.\ \ref{T:approximate} depends upon $p$.  In particular, the proof we shall give does not extend to $p=\infty$. 
This leaves open the following question:
Can one find an $L^\infty$-$\delta$-approximate critical hermitian structure in general?
\end{Rem}

The proof of Thm.\ \ref{T:approximate} will be given in Section \ref{S:proofapproximate} below.    A preliminary result in this direction is obtained by passing to a resolution of the filtration.  We will prove the following:
\begin{Prop} \label{P:approximate}
Let $E$, $X$, and $\omega$ be as in Thm.\ \ref{T:approximate}.  Let 
$
\mu_i=\mu_\omega\left(\filt_i(E)/\filt_{i-1}(E)\right)
$.
Then there is a sequence of monoidal transformations giving a K\"ahler surface $\pi: \widehat X\to X$, a number $p_0>1$, and a
family of K\"ahler metrics $\omega_\varepsilon$ converging to $\pi^\ast \omega$ as $\varepsilon\to 0$, such that the following holds:  Let
$\widehat {\mathcal F}$ be the filtration of $\pi^\ast E=\widehat E$ given by $\left\{  \sat_{\widehat E}(\pi^\ast\filt_i(E))\right\} $.  Then for
any $\delta >0$ and any $1\leq p<p_0$ there is $\varepsilon_1>0$ and a smooth hermitian metric $\widehat H$ on $\widehat E$  such that
for all $0<\varepsilon\leq \varepsilon_1$,
$
\Vert\sqrt{-1}\, \Lambda_{\omega_\varepsilon} F_{(\dbar_{\widehat E}, \widehat H)}-\Psi(\widehat
{\mathcal F}, (\mu_1,\ldots,\mu_\ell), \widehat H)\Vert_{L^p(\omega_\varepsilon)}\leq
\delta
$.
\end{Prop}

 For a \emph{fixed} $\varepsilon>0$ sufficiently small compared to $\delta$, the analogous result for \emph{any} $p$ is a consequence of Thm.\
\ref{T:bundlecase}.  The key point in
the statement of Prop.\ \ref{P:approximate} is that a metric $\widehat H$ may be found which satisfies the condition of the proposition
\emph{uniformly} in $\varepsilon$.  
 The direct construction of $\widehat H$ given below, however,  requires that $p$ be sufficiently small.  This requirement derives
from the following:

\begin{Lem} \label{L:approximate}
Let $\pi : \widehat X\to X$ be a blow-up of the type discussed in Section \ref{S:hnresolution}, and let
$\omega_\varepsilon=\pi^\ast\omega+\varepsilon\eta$ be the family of K\"ahler metrics defined in Lemma \ref{L:kahler}.  Then there is associated to
$\widehat X$ a positive integer $\hat m$ with the following property:  given any
$p$, 
$
 1\leq p< 1+ (1/\hat m)
$,
 there is $\varepsilon_1>0$ such that for any $\tilde p$ satisfying
$ p(1-\hat m(p-1))^{-1}< \tilde p \leq +\infty 
$,
there is a constant $C(\tilde p,\varepsilon_1)$ such that
$
\left\Vert \Lambda_{\omega_\varepsilon} G \right\Vert_{L^p(\omega_\varepsilon)}\leq C(\tilde p, \varepsilon_1)
\left\Vert \Lambda_{\omega_{\varepsilon_1}} G \right\Vert_{L^{\tilde p}(\omega_{\varepsilon_1})}
$,
for all smooth $(1,1)$ forms $G$ on $\widehat X$ and all $0<\varepsilon\leq \varepsilon_1$.
\end{Lem}

\begin{proof}  Since $\omega_\varepsilon\to\pi^\ast\omega$ smoothly, and $\pi^\ast\omega$ is a K\"ahler metric off the exceptional set, the estimate
is clearly local near the exceptional set {\bf e}.  Let $\hat x\in{\bf e}\subset\widehat X$ with $x=\pi(\hat x)$.  We may choose local coordinates
$(z_1,z_2)$ near $x$ with respect to which the K\"ahler form $\omega$ is standard to first order.
Since {\bf e} has normal crossings, we may choose coordinates $(\xi_1,\xi_2)$ in a neighborhood $\widehat U$ of $\hat x$, centered at $\hat x$, and 
such that ${\bf e}\cap\widehat U$ is contained in the union of the coordinate axes $\{\xi_1=0\}\cup\{\xi_2=0\}$.  Regarding $z_1, z_2$ as holomorphic
functions on $\widehat U$, let us write:
$z_1 \simeq \xi_1^a\xi_2^m $,
$z_2 \simeq \xi_1^b\xi_2^n$, modulo higher order terms, where
 $a,b,m,n$ are nonnegative integers and $an\neq bm$.
In these coordinates we have:
\begin{equation}  \label{E:pullback}
\pi^\ast(\omega\wedge\omega)\simeq|\xi_1|^{2(a+b-1)}|\xi_2|^{2(m+n-1)}(an-bm)^2\left(\tfrac{\sqrt{-1}\, }{2} d\xi_1\wedge d\bar \xi_1\right)
\wedge \left(\tfrac{\sqrt{-1}\, }{2} d\xi_2\wedge d\bar \xi_2\right)
 \end{equation}
 modulo  higher order terms.
 At this point we set:
$
\hat m=\max\left\{ (a+b-1), (m+n-1)\right\}
$.
Let $g^\varepsilon_{\alpha\bar\beta}$ denote the K\"ahler metric with K\"ahler form $\omega_\varepsilon$ in the coordinates $\xi_\alpha$, $\alpha=1,2$ on the neighborhood $\widehat U$.  It
follows from (\ref{E:pullback})  that there is a constant $C$ uniform in $\varepsilon$ such that: 
\begin{equation} \label{E:determinant}
|\det g^\varepsilon_{\alpha\bar\beta} | \geq C|\xi_1|^{2\hat m}|\xi_2|^{2\hat m}\ .
\end{equation}
Let $g_\varepsilon^{\alpha\bar\beta}=M_\varepsilon^{\alpha\bar\beta}/\det g^\varepsilon_{\alpha\bar\beta}$ denote the inverse metric.  Since
$\dim_{\CBbb} X=2$, $M_\varepsilon^{\alpha\bar\beta}=\epsilon^{\alpha\sigma}\epsilon^{\beta\rho}g^\varepsilon_{\sigma\bar\rho}$, where
$\epsilon^{11}=\epsilon^{22}=0$, $\epsilon^{12}=-\epsilon^{21}=1$.  In particular, we can arrange for a constant $C$ uniform in $\varepsilon$ such
that in a neighborhood of $\hat x$ (still denoted $\widehat U$):
\begin{equation} \label{E:M}
|M_\varepsilon^{\alpha\bar\beta}|\leq C|M_{\varepsilon_1}^{\alpha\bar\beta}|\ ,
\end{equation}
for all $0<\varepsilon\leq \varepsilon_1$.  
Using (\ref{E:determinant}) and (\ref{E:M}) we now prove the lemma.  Let $G_{\alpha\bar\beta}$ be the local expression in the coordinates
$(\xi_1,\xi_2)$ on
$\widehat U$ of the form
$G$ in the statement of the lemma.  Then by (\ref{E:M}) there is a constant $C$ independent of $\varepsilon$ and $G$ such that:
$
|G_{\alpha\bar\beta}M_\varepsilon^{\alpha\bar\beta}|^p\leq C|G_{\alpha\bar\beta}M_{\varepsilon_1}^{\alpha\bar\beta}|^p
$,
for all $p\geq 1$ and $0<\varepsilon\leq \varepsilon_1$.  Then:
\begin{align*}
\left\Vert \Lambda_{\omega_\varepsilon} G \right\Vert^p_{L^p(\widehat U, \omega_\varepsilon)}
&=\int_{\widehat U}|G_{\alpha\bar\beta}g_\varepsilon^{\alpha\bar\beta}|^p \left(\det g^\varepsilon_{\alpha\bar\beta}\right) |d\xi_1|^2|d\xi_2|^2 \\
&\leq C\int_{\widehat U}|G_{\alpha\bar\beta}g_{\varepsilon_1}^{\alpha\bar\beta}|^p \left(\frac{\det
g^\varepsilon_{\alpha\bar\beta}}{\det
g^{\varepsilon_1}_{\alpha\bar\beta}}\right)^{1-p}\left(\det
g^{\varepsilon_1}_{\alpha\bar\beta}\right)|d\xi_1|^2|d\xi_2|^2
\end{align*}
Now $\det g^{\varepsilon_1}_{\alpha\bar\beta}$ is uniformly bounded away from zero on $\widehat U$ by a constant depending upon $\varepsilon_1$.  If $\tilde p=\infty$, the result follows from  
 (\ref{E:determinant}) and the assumption on $p$.  If $\tilde p\neq\infty$, we apply H\"older's inequality with the conjugate variables: $r=\tilde p/p$, $s=\tilde p/(\tilde p-p)$, and find:
 \begin{equation} \label{E:int}
 \Vert \Lambda_{\omega_\varepsilon} G \Vert^p_{L^p(\widehat U, \omega_\varepsilon)}\leq C\Vert \Lambda_{\omega_{\varepsilon_1}} G
\Vert^p_{L^{\tilde p}(\widehat U, \omega_{\varepsilon_1})}
\biggl\{\int_{\widehat U}(\det g^{\varepsilon}_{\alpha\bar\beta})^{(1-p)s}
|d\xi_1|^2|d\xi_2|^2\biggr\}^{1/s} 
 \end{equation}
From the assumption on $\tilde p$ we have $2\hat m (1-p)s>-2$.  By
(\ref{E:determinant}), this implies that the integral on the right hand side of (\ref{E:int}) is convergent uniformly in $\varepsilon$.  This proves the lemma.
\end{proof}

\begin{proof}[Proof of Prop.\ \ref{P:approximate}]
We proceed by induction on the rank of $E$.  The case of rank 1 is trivial, since line bundles admit Hermitian-Einstein metrics.  Suppose
that $\rk(E) > 1$, and consider the HN filtration $\{\filt_{i,\omega}(E)\}$.  If the filtration is by subbundles, then the result follows
from Thm.\ \ref{T:bundlecase}.  Consider the case where the filtration is not by subbundles.  For convenience, set
$E_i=\filt_{i,\omega}(E)$, $Q_i=E_i/E_{i-1}$, and $\mu_i=\mu_\omega(Q_i)$.  By  Prop.\  \ref{P:resolution2} there is a resolution $\pi:\widehat X\to
X$ where the filtration
$\widehat E_i=\sat_{\widehat E}\pi^\ast E$ is a filtration of $\widehat E=\pi^\ast E$ by subbundles.  Let $\widehat Q_i=\widehat E_i/\widehat E_{i-1}$.
By the inductive hypothesis, given $\delta>0$ and for any $\varepsilon >0$ sufficiently small we may find $L^{\hat p}$ $\delta$-approximate
critical hermitian structures
$\widehat H^\varepsilon_i$ on each $\widehat Q_i$, for some $\hat p>1$. Since $E_i/E_{i-1}$ is semistable, it follows from Prop.\ \ref{P:stability} that for a given $\delta_1$ we may assume $\varepsilon_1$ has been
chosen such that
$
| \mu_{max,\varepsilon}(\widehat Q_i)-\mu_{min,\varepsilon}(\widehat Q_i)| \leq \delta_1
$
for all $0<\varepsilon\leq \varepsilon_1$. Here, $\delta_1>0$ will be chosen presently. In particular:
\begin{equation}  \label{E:phibound}
\left\Vert \Psi^{hn}_{ \omega_\varepsilon}(\dbar_{\widehat Q_i}, \widehat H^{\varepsilon_1}_i) -\mu_i{\bf I}_{\widehat Q_i}\right\Vert_{L^{\hat p}(\omega_\varepsilon)}
\leq C\delta_1\ ,
\end{equation}
for a constant $C$ independent of $\varepsilon$ and $\delta_1$.
Associated to $\widehat X$ is an integer $\hat m$ as in Lemma \ref{L:approximate}.  We choose $p_0$ sufficiently close to $1$ so
that
$
p_0 < \hat p (1-\hat m(p_0-1))
$.
Then the conclusion of the lemma, along with (\ref{E:phibound}), guarantee that for each $1\leq p<p_0$, each $i$, and each
$0<\varepsilon\leq
\varepsilon_1$,
$
\Vert\sqrt{-1}\, \Lambda_{\omega_\varepsilon} F_{(\dbar_{\widehat Q_i}, \widehat H_i^{\varepsilon_1})}-\mu_i{\bf I}_{\widehat
Q_i}\Vert_{L^p(\omega_\varepsilon)}\leq
C\delta_1
$,
for a constant $C$ independent of $\varepsilon$ and $\delta_1$.  Choose a smooth splitting $\widehat E=\oplus_{i=1}^\ell \widehat Q_i$, and 
let $\widehat H=\oplus_{i=1}^\ell \widehat H_i^{\varepsilon_1}$.  This is a smooth metric on $\widehat E$.  Since the filtration $\{\widehat E_i\}$ is
by subbundles, we may argue as in the proof of Thm.\ \ref{T:bundlecase} that $\widehat H$ may be modified to produce the desired result if we
choose $\delta_1$ sufficiently small (depending upon the constant $C$) compared to $\delta$.
\end{proof}

%%%%%%%%%%%%%%%%%%%%%%%%%%%%%%%%%%%%%%%%%%%%%%%%%%%%%%%%%%%%%%%%%%%%%%%%%%%%%%%%%%%%%%%%%%%%%%%%%%%%%%%%%%%%%%%%%%
%%%%%%%%%%%%%%%%%%%%%%%%%%%%%%%%%%%%%%%%%%%%%%%%%%%%%%%%%%%%%%%%%%%%%%%%%%%%%%%%%%%%%%%%%%%%%%%%%%%%%%%%%%%%%%%%%%

\section{The Harder-Narasimhan Type of the Uhlenbeck Limit}       \label{S:type}

%%%%%%%%%%%%%%%%%%%%%%%%%%%%%%%%%%%%%%%%%%%%%%%%%%%%%%%%%%%%%%%%%%%%%%%%%%%%%%%%%%%%%%%%%%%%%%%%%%%%%%%%%%%%%%%%%%
%%%%%%%%%%%%%%%%%%%%%%%%%%%%%%%%%%%%%%%%%%%%%%%%%%%%%%%%%%%%%%%%%%%%%%%%%%%%%%%%%%%%%%%%%%%%%%%%%%%%%%%%%%%%%%%%%%

As indicated in the Introduction, the proof of the Main Theorem proceeds in two steps.  The goal of this section is to prove the first step:

\begin{Thm} \label{T:type}
Let $D_{t}$ be  a solution to the YM flow equations (\ref{E:ymflow}) with initial condition $D_0$ and  Uhlenbeck limit $D_\infty$.  Let $E_\infty$
denote the holomorphic vector bundle obtained from $D_\infty$, as in Prop.\ \ref{P:convergence}.  Then the Harder-Narasimhan type of $(E_\infty,
D_{\infty}^{\prime\prime})$ is the same as that of $(E,D_0^{\prime\prime})$.
\end{Thm}

\noindent We prove this theorem in the first subsection below.  In the second subsection, we use this fact to prove Thm.\ \ref{T:approximate}.
 
%%%%%%%%%%%%%%%%%%%%%%%%%%%%%%%%%%%%%%%%%%%%%%%%%%%%%%%%%%%%%%%%%%%%%%%%%%%%%%%%%%%%%%%%%%%%%%%%%%%%%%%%%%%%%%%%%%%

\subsection{Proof of Thm.\ \ref{T:type}}   \label{S:hartman}

We begin with the following:
    
\begin{Lem} \label{L:close}
Let $E\to X$ be a holomorphic bundle of HN type $\vec \mu_0$. There is $\alpha_0>1$ such that the following holds:   given any $\delta>0$ and any $N$,  there is a
hermitian metric
$H$ on $E$ such that
$\HYM_{\alpha, N}(\dbar_E,H)\leq
\HYM_{\alpha, N}(\vec\mu_0)+\delta$,
 for all $1\leq \alpha\leq
\alpha_0$.
\end{Lem}
   
   \begin{proof}
   To begin,
let $\pi:\widehat X\to X$ be a resolution  of the HN filtration guaranteed by Prop.\ \ref{P:resolution2}, $\omega_\varepsilon$ the family of K\"ahler metrics from Lemma \ref{L:kahler},
    and  $\widehat
E=\pi^\ast E$.  As a  direct consequence Prop.\ \ref{P:approximate}, where $\hat\alpha_0=p_0$, and Cor.\ \ref{C:hntype}, there is $\hat\alpha_0>1$ such that the following holds:  given any $\delta>0$ there exists a smooth hermitian metric $\widehat H$ on
$\widehat E$,  and $\varepsilon_1>0$  (depending on $\widehat H$) such that:
\begin{equation} \label{E:initialestimate}
\HYM^{\omega_\varepsilon}_{\alpha, N}(\dbar_{\widehat E},\widehat H)\leq
\HYM_{\alpha, N}(\vec\mu_0)+\delta/2\ ,
\end{equation}
for all $1\leq \alpha\leq
\hat\alpha_0$, and all $0\leq \varepsilon\leq\varepsilon_1$.
In order to obtain a metric on $X$, we use a cut-off argument.  Let $x\in \algsing$, and choose a coordinate neighborhood $U$ of $x$.
    For $R>0$ sufficiently small, let $B_R\subset U$ denote the ball of radius $R$ about $x$ with respect to these coordinates.   We also choose a holomorphic  trivialization of $E\to U$.   This gives a trivialization of $\widehat E$ on $\widehat U=\pi^{-1}(U)$, with respect to which we regard $\widehat H$ as a positive definite hermitian matrix valued function.
    Given $R$, 
     we may choose a smooth function $\varphi_R$ on $U$, $0\leq \varphi_R\leq 1$, $\varphi_R\equiv 0$ on a ball of radius $R/2$ centered at $p$, and $\varphi_R\equiv 1$ outside a ball of radius $R$, and such that $|\varphi_R'|\leq CR^{-1}$ and $|\varphi_R^{\prime\prime}|\leq CR^{-2}$, where $C$ is a constant independent of $R$.  Define a metric $H_{\varphi_R}$ as follows:  If $\widehat H s_i=\lambda_i s_i$ with respect to a unitary frame $\{s_i\}$, then $H_{\varphi_R} s_i=(\varphi_R\lambda_i+(1-\varphi_R))s_i$.
 With this definition, $H_{\varphi_R}$ extends as a smooth metric on $E\to U$.  Let $\widehat   H_{\varphi_R}$ denote the pull-back metric on $\widehat X$.  A calculation then shows that there are constants $C_1$ and $C_2$, depending on $\widehat H$ but not on $R$ or $\varepsilon$, such that on $\pi^{-1}(B_R\setminus B_{R/2})$:
 $
  \Lambda_{\omega_\varepsilon} F_{(\dbar_{\widehat E}, \widehat   H_{\varphi_R})}=
 \Lambda_{\omega_\varepsilon} F_{(\dbar_{\widehat E},\widehat H)} + f_{\varepsilon, R}
 $,
 where in the  coordinates used in the proof of Lemma \ref{L:approximate}:
 \begin{equation}\label{E:fbound}
 |f_{\varepsilon,R}|\leq |\det g^\varepsilon_{\alpha\bar\beta}|^{-1}(C_1+C_2 R^{-2})\ .
 \end{equation}
 Continuing  this way for all points in $\algsing$, we obtain a metric, still denoted $H_{\varphi_R}$,   with $\widehat H_{\varphi_R}=\widehat H$ outside  the union $U_R$ of the balls $B_R$,  $\vol_\omega(U_R)\simeq R^4$, and $H_{\varphi_R}$ standard with respect to the trivialization inside $U_{R/2}$.  Hence:
   \begin{align*}
| \HYM^{\omega_\varepsilon}_{\alpha, N}(\dbar_{\widehat E},\widehat H_{\varphi_R})&- \HYM^{\omega_\varepsilon}_{\alpha, N}(\dbar_{\widehat E},\widehat H)|
 \leq C R^4 \\
&\qquad \qquad+ C\Vert \Lambda_{\omega_\varepsilon} F_{(\dbar_{\widehat E},\widehat H)}\Vert^{\alpha}_{L^\alpha_{\omega_\varepsilon}(\pi^{-1}(U_R))} 
+ C\Vert f_{\varepsilon,R}\Vert^{\alpha}_{L^\alpha_{\omega_\varepsilon}(\pi^{-1}(U_R\setminus U_{R/2}))}
\end{align*}
where $C$ is independent of $R$ and $\varepsilon$.   By the construction of $\widehat H$, the second term on the right hand side tends to zero as $R\to 0$, uniformly in $\varepsilon$.  Hence, we may choose $R$ sufficiently small so that this term is less than $\delta/4$, say, for all $\varepsilon\leq \varepsilon_1$.
Letting $\varepsilon\to 0$ and using \eqref{E:fbound} to bound the third term,  we obtain  an estimate of the form:
$$
| \HYM^{\omega}_{\alpha, N}(\dbar_{E}, H_{\varphi_R})- \HYM^{\omega_0}_{\alpha, N}(\dbar_{\widehat E},\widehat H)|\leq  C(1+R^{-2\alpha})R^4+\delta/4\ .
$$
Now by \eqref{E:initialestimate}, provided $\alpha_0\leq\hat\alpha_0$ and $\alpha_0<2$, we may take  $R$ sufficiently small so that:
$$
\HYM^{\omega}_{\alpha, N}(\dbar_{ E}, H_{\varphi_R})
    \leq
\HYM_{\alpha, N}(\vec\mu_0)+\delta
   $$
 for all $\alpha\leq \alpha_0$. 
 \end{proof}
 
 Next, we have the following ``distance decreasing" result:
 
 \begin{Lem} \label{L:hymconvergence}
Let $\alpha_0$ be as in Lemma \ref{L:close}.   Let $H$ be any smooth hermitian metric on $E$, and let $D_t$ be a solution to
the YM flow equations (\ref{E:ymflow}) with initial condition $(\dbar_E, H)$.  Then:  $$\lim_{t\to\infty}\HYM_{\alpha, N}(D_t)=\HYM_{\alpha,N}(\vec\mu_0)\ ,$$
  for all $1\leq \alpha\leq\alpha_0$, and all $N$.  In particular, if $D_\infty$ is the
Uhlenbeck limit along the flow, then 
$\HYM_{\alpha,N}(D_\infty)=\HYM_{\alpha,N}(\vec\mu_0)$.
\end{Lem}
 
 \begin{proof}
 We first point out that the second assertion follows from the first because of Prop.\ \ref{P:equality}.
  For fixed $\alpha$, $1\leq\alpha\leq \alpha_0$, and fixed $N$, define $\delta_0>0$ by:
\begin{equation} \label{E:delta0}
2\delta_0+\HYM_{\alpha,N}(\vec\mu_0)=\min \{\HYM_{\alpha,N}(\vec\mu) : \HYM_{\alpha,N}(\vec\mu)>\HYM_{\alpha,N}(\vec\mu_0)\}\ ,
\end{equation}
where $\vec\mu$ runs over all possible HN types of holomorphic vector bundles on $X$ with the rank of $E$.  
Consider  metrics $H$ on $E$ with associated connection $D=(\dbar_E, H)$ satisfying:
  \begin{equation} \label{E:startclose}
\HYM_{\alpha,N}(D)\leq\HYM_{\alpha,N}(\vec\mu_0)+\delta_0\ .
\end{equation}
  Let $D_{\infty}$ be the Uhlenbeck limit along the flow with
initial condition $D$.  Then combining Prop.\ \ref{P:slopesemicontinuity}, Prop.\ \ref{P:key1} (1), and Prop.\ \ref{P:key2}, we have:
$$
\HYM_{\alpha,N}(\vec\mu_0)\leq \HYM_{\alpha,N}(D_\infty)\leq 
\HYM_{\alpha,N}(D)\leq\HYM_{\alpha,N}(\vec\mu_0)+\delta_0\ .
$$
Hence, by \eqref{E:delta0} we must have $\HYM_{\alpha,N}(D_\infty)=\HYM_{\alpha,N}(\vec\mu_0)$.   This shows that the result holds for initial conditions satisfying \eqref{E:startclose}.     
    
    In the following, let us denote by $ D^H_t$ the solution to the YM flow at time $t$  with initial condition $
    D^H_0=(\dbar_E, H)$.
We are going to prove that for \emph{any} initial condition $H$ and any $\delta>0$, there is $T\geq 0$ such that:
    \begin{equation} \label{E:hyminequality}
\HYM_{\alpha,N}(D^H_t) < \HYM_{\alpha,N}(\vec\mu_0)+\delta\ ,\qquad\text{for all}\ t\geq T\ .
\end{equation}
  Without loss of generality, assume $0<\delta\leq\delta_0/2$.  Let
$
{\mathcal H}_\delta
$ denote the set of smooth hermitian metrics  $ H$ on $ E$ with the  property that (\ref{E:hyminequality}) holds for $D_t^{H}$ and some $T$. From the discussion above, ${\mathcal H}_\delta$ is nonempty: indeed, any metric satisfying \eqref{E:startclose} is in ${\mathcal H}_\delta$, and according to Lemma \ref{L:close} we may always find such a metric.
% Let
%$
%M_\delta=\sup\left\{\HYM_{\alpha,N}( \dbar_E,H) :  H\in {\mathcal H}_\delta\right\}
%$.
% We need to show that $M_\delta=+\infty$.  Suppose to the contrary that $M_\delta<+\infty$ for some
%$\delta$. Let us call a metric  ``bad" if there is no $T$ for which (\ref{E:hyminequality}) holds,  and ``good" otherwise. 
%Choose sequences 
 Let $H^j$ be  a sequence of smooth hermitian metrics on $ E$ such that each $ H^j\in {\mathcal H}_\delta$,  and 
suppose $ H^j\to  K$,  in the $ C^\infty$ topology, for some metric $K$.  
    Since $ H^j\in{\mathcal H}_\delta$ we have a sequence $T_j$ such that for  all $t\geq T_j$:
\begin{equation} \label{E:good}
\HYM_{\alpha,N}( D^{ H^j}_t)\leq\HYM_{\alpha,N}( D^{ H^j}_{T_j})\leq
\HYM_{\alpha,N}(\vec\mu_0)+\delta\ .
\end{equation}
     By Lemma \ref{L:subsolution} (2) and the $C^\infty$ convergence of $H^j$, $\Vert\Lambda F_{D^{ H^j}_{t}}\Vert_{L^\infty}$ and $\Vert\Lambda F_{D^{ K}_{t}}\Vert_{L^\infty}$ are bounded  uniformly  for all $t\geq 0$. Hence,   it follows from Prop.\ \ref{P:ymcompactness}  that we may find a sequence $t_j\geq T_j$, Yang-Mills connections $D^{(1)}_\infty$ and $D^{(2)}_\infty$, and bubbling sets $\ansing_{(1)}$
     and $\ansing_{(1)}$, such that 
    $D^{ H^j}_{t_j}\weakarrow D^{(1)}_\infty$ in $L^p_{1,loc.}(X\setminus \ansing_{(1)})$ and
    $D^{ K}_{t_j}\weakarrow D^{(2)}_\infty$ in $L^p_{1,loc.}(X\setminus \ansing_{(2)})$, for all $1\leq p<\infty$.  Moreover, by Cor.\ \ref{C:lpconvergence},  $\Lambda F_{ D^{ H^j}_{t_j}}\to
\Lambda F_{ D^{(1)}_\infty}$  and $\Lambda F_{ D^{ K}_{t_j}}\to
\Lambda F_{ D^{(2)}_\infty}$ strongly in $L^p$, for all $p$.
    
    \begin{Claim}  
   $ D^{(1)}_\infty=D^{(2)}_\infty$.
    \end{Claim}
    
\begin{proof}[Proof of the Claim]
Write $ H^j_{t_j}= h^j_{t_j} K_{t_j}$.  It follows by \cite[Prop.\ 13]{Do1} that
$
 \sup\sigma( H^j_t,  K_t)\to 0
$
as $j\to \infty$, uniformly in $t$, where:
    $
    \sigma(H,K)=\tr H^{-1}K +\tr K^{-1}H -2\rk(E)
    $,
is the usual $C^0$-distance on the space of hermitian metrics on $E$.  In particular,
$\sup| h^j_{t_j}-{\bf I}_E|\to 0$ as $j\to
\infty$.
Let
$\ansing=\ansing_{(1)}\cup \ansing_{(2)}$, and choose a smooth test form $\phi\in\Omega^{1,0}(\End E)$, compactly supported on $X\setminus \ansing$. 
 We have 
$
(D^{ H^j}_{t_j})'- (D^{ K}_{t_j})'=( h^j_{t_j})^{-1} ( D^{ K}_{t_j})^\prime
( h^j_{t_j})
$.
For notational simplicity, set
$ D_j= D^{ K}_{t_j}$, and $ h_j=h^j_{t_j}$.  Then there is a constant $C$ such that
$$
\bigl|\langle  h_j^{-1}  D_j^\prime ( h_j), \phi\rangle_{L^2}\bigr| \leq C\bigl|\langle  h_j, ( D_j^\prime)^\ast 
\phi\rangle_{L^2}\bigr| 
   \leq C\left\{
\bigl|\langle h_j, ( D_j^\prime- D_{\infty}^\prime)^\ast 
\phi\rangle_{L^2}\bigr|
+
\bigl|\langle h_j, ( D_{\infty}^\prime)^\ast 
\phi\rangle_{L^2}\bigr|
\right\}\ .
$$
Now $ D_j^\prime \weakarrow  D_{\infty}^\prime$ in $L^p_{1, loc.}$, so we may assume 
    $ D_j^\prime \to  D_{\infty}^\prime$  in $C^0$. 
Combined with the uniform bound for $\Vert h_j\Vert_{L^\infty}$, this implies that the first term  on the right hand side above goes to
zero.  For the second term, notice that since $ h_j\overset{\scriptscriptstyle C^0}{\lra}{\bf I}_E$:
$$
\langle  h_j, ( D_{\infty}^\prime)^\ast 
\phi\rangle_{L^2}\lra \langle {\bf I}_E, ( D_{\infty}^\prime)^\ast 
\phi\rangle_{L^2}
=\int_X \tr ( D_{\infty}^\prime)^\ast 
\phi\, dvol
=\int_X \partial^\ast \tr 
\phi\, dvol= 0\ ,
$$
by Stokes' theorem.  This proves that $(D^{ H^j}_{t_j}- D^{ K}_{t_j})\weakarrow 0$ in $L^2_{loc.}(X\setminus\ansing)$, and  the claim follows.
    \end{proof}

Set  $ D_\infty=D^{(1)}_\infty =  D^{(2)}_\infty$. 
Since $\Lambda F_{ D^{ H^j}_{t_j}}\to
\Lambda F_{ D_\infty}$  and $\Lambda F_{ D^{ K}_{t_j}}\to
\Lambda F_{ D_\infty}$ strongly in $L^p$, for all $1\leq p<\infty$, we have (see Lemma \ref{L:norm}):
$$
\lim_{j\to\infty}\HYM_{\alpha,N}( D^{ H^j}_{t_j})=
\lim_{j\to\infty}\HYM_{\alpha,N}( D^{ K}_{t_j})=
\HYM_{\alpha,N}( D_\infty)\ .
$$
Hence, for $j$ sufficiently large:
\begin{align*}
\HYM_{\alpha,N}( D^{ K}_{t_j})
&\leq
\HYM_{\alpha,N}( D_\infty) +\delta 
=
\lim_{j\to\infty} \HYM_{\alpha}( D^{ H^j}_{t_j})+\delta \\
&\leq
\HYM_{\alpha,N}(\vec\mu_0) +2\delta
\leq
\HYM_{\alpha,N}(\vec\mu_0) +\delta_0 \ ,
    \end{align*}
    where in the second line we have used \eqref{E:good} and the fact that
 $ \delta\leq\delta_0/2$.
It follows from \eqref{E:delta0} and the discussion above that for $j$ sufficiently large:
$\lim_{t\to\infty}\HYM_{\alpha,N}( D^{ K}_{t_j+t})= \HYM_{\alpha,N}\left(\vec\mu_0\right)$.  In particular,
$
\HYM_{\alpha,N}( D^{ K}_{t_j+t})< \HYM_{\alpha,N}\left(\vec\mu_0\right)+\delta
$, for $t\geq T$, $T$ sufficiently large.
Therefore, $ K\in {\mathcal H}_\delta$.  This proves that ${\mathcal H}_\delta$ is closed in the $C^\infty$ topology.  By the continuous dependence of the flow on initial conditions, ${\mathcal H}_\delta$ is also open. Since the space of smooth metrics is connected, we conclude that every metric is in ${\mathcal H}_\delta$, 
and \eqref{E:hyminequality} holds for all $\delta>0$ and all initial conditions $H$.  In particular, we can choose $\delta\leq \delta_0$ and conclude that  $\lim_{t\to\infty}\HYM_{\alpha,N}( D^{H}_{t})= \HYM_{\alpha,N}\left(\vec\mu_0\right)$, for any $H$.
Since the choice of $N$ was arbitrary, the proof is complete.
\end{proof}

Finally, we have:

\begin{proof}[Proof of Thm.\ \ref{T:type}]
Let $\vec\mu_0=(\mu_1,\ldots,\mu_R)$ (resp.\ $\vec\lambda_\infty=(\lambda_1,\ldots,\lambda_R)$) be the HN type of $(E,D_0^{\prime\prime})$ (resp.\ $(E,D_\infty^{\prime\prime})$).  
By Lemma \ref{L:hymconvergence}, $\varphi_\alpha(\vec\mu_0+N)=\varphi_\alpha(\vec\lambda_\infty+N)$ for all $1\leq \alpha\leq \alpha_0$ and all $N$.  In particular, we may
choose $N$ sufficiently large so that $\mu_R+N\geq 0$.  By Prop.\ \ref{P:slopesemicontinuity} we also have $\lambda_R+N\geq 0$.  Since $\alpha_0>0$, the hypotheses of Prop.\ \ref{P:key1} (2) are then satisfied, and we conclude that $\vec\mu_0+N=\vec\lambda_\infty+N$, and so $\vec\mu_0=\vec\lambda_\infty$.
\end{proof}

%%%%%%%%%%%%%%%%%%%%%%%%%%%%%%%%%%%%%%%%%%%%%%%%%%%%%%%%%%%%%%%%%%%%%%%%%%%%%%%%%%%%%%%%%%%%%%%%%%%%%%%%%%%%%%%%%%

\subsection{Proof of Theorem \ref{T:approximate}}  \label{S:proofapproximate}
   
   Let $(E, D_0^{\prime\prime})$ be a holomorphic bundle, $D_j=g_j(D_0)$ a sequence of unitary connections in the $\GC$ orbit of $D_0$, and set $F_j=F_{D_j}$.  For the next result we make the following assumptions:
        \begin{Assume}  \label{A:1}
        \begin{enumerate}
     \item $D_j$ converges off a finite set of points $\ansing\subset X$ weakly in $L^p_{1,loc.}$, for all $p>4$, to a Yang-Mills connection $D_\infty$ on a bundle $E_\infty$;
  \label{convergence}
   \item The HN type of $(E_\infty, D_\infty^{\prime\prime})$ is the same as the HN type of $(E, D_0^{\prime\prime})$; \label{type}
       \item  $\Vert \Lambda F_j\Vert_{L^\infty}$ is bounded uniformly in $j$, and $\Lambda F_{j}\overset{\scriptscriptstyle L^1}{\lra} \Lambda F_{\infty} $, where $F_\infty=F_{D_\infty}$.  \label{bound}
   \end{enumerate}
   \end{Assume}
   
   Recall that for a weakly holomorphic projection $\pi$ of $E$, the rank and degree of $\pi$ are, by definition, the rank and degree of the associated saturated  subsheaf of $E$ (see the discussion in the proof of Lemma \ref{L:degreesemicontinuity}).
      
   \begin{Lem} \label{L:projectionconvergence}
 \begin{enumerate}
 \item  Let $\{\pi_j^{(i)}\}$ be the HN filtration of $(E,D_j^{\prime\prime})$
 and $\{\pi_\infty^{(i)}\}$ the HN filtration of $(E_\infty,D_\infty^{\prime\prime})$.  Then after passing to a subsequence,
 $\pi_j^{(i)}\to \pi_\infty^{(i)}$ strongly in $L^p\cap L^2_{1,loc.}$, for all $1\leq p<\infty$ and all $i$.
 \item  Suppose $(E,D_0^{\prime\prime})$ is semistable and $\{\pi_{ss,j}^{(i)}\}$ are Seshadri filtrations of $(E,D_j^{\prime\prime})$.  Without loss of generality, assume the ranks of the subsheaves $\pi_{ss,j}^{(i)}$ are constant in $j$.  Then there is a filtration $\{\pi_{ss,\infty}^{(i)}\}$ of $(E_\infty,D_\infty^{\prime\prime})$ such that
 after passing to a subsequence,
 $\pi_{ss,j}^{(i)}\to \pi_{ss,\infty}^{(i)}$ strongly in $L^p\cap L^2_{1,loc.}$, for all $1\leq p<\infty$ and all $i$.  Moreover, the rank and degree of  $\pi_{ss,\infty}^{(i)}$ is equal to the rank and degree of $\pi_{ss,j}^{(i)}$ for all $i$ and $j$.
 \end{enumerate}
 \end{Lem}
 
 \begin{proof}   For part (1), set $E_i=\filt_i(E, D_0^{\prime\prime})$ and $E_\infty^{(i)}=\filt_i(E_\infty, D_\infty^{\prime\prime})$.  Hence, $\pi_j^{(i)}$ is the orthogonal
projection onto the subsheaf $g_j(E_i)$ of $(E,D_j^{\prime\prime})$.  As in the proof of Lemma \ref{L:degreesemicontinuity} we have:
$$
\deg(E_i)+\frac{1}{2\pi}\int_X\Vert D_j^{\prime\prime}\pi_j^{(i)}\Vert^2\, dvol\leq
\sum_{k\leq \rk(E_i)}\mu_k+\frac{1}{2\pi}\Vert\Lambda F_{j}-\Lambda F_{\infty}\Vert_{L^1}\ ,
$$
 where $\vec\mu=(\mu_1,\ldots,\mu_R)$ is the  HN type of $(E_\infty,D_\infty^{\prime\prime})$. By the assumption \P\  \ref{A:1} \eqref{type}, $\vec\mu$ is also the HN type of $(E,D_0^{\prime\prime})$, so
$
\deg(E_i)=\sum_{k\leq \rk(E_i)}\mu_k
$.
It then follows from \P\ \ref{A:1} \eqref{bound} that:
 \begin{equation} \label{E:ffconvergence}   
    D_j^{\prime\prime}\pi_j^{(i)}\to 0\ \text{in}\ L^2\ .
    \end{equation}
    Write
     \begin{equation} \label{E:bridge}
 D_\infty^{\prime\prime}\pi_j^{(i)}=D_j^{\prime\prime}\pi_j^{(i)}+( D_\infty^{\prime\prime}-D_j^{\prime\prime})\pi_j^{(i)}\ .
\end{equation}
We may assume,  after perhaps 
     passing to a
subsequence, that $\pi_j^{(i)}\weakarrow \tilde\pi_\infty^{(i)}$ in $L^2_{1,loc.}$, for some $L^2_1$ projection $\tilde \pi_\infty^{(i)}$.  
 Since $\pi_j^{(i)}$ is uniformly bounded, $\pi_j^{(i)}\overset{\scriptscriptstyle L^p}{\lra}\tilde \pi_\infty^{(i)}$ for all $p$.  Then  as in the proof of Prop.\ \ref{P:ymcompactness} we conclude from (\ref{E:ffconvergence}) and (\ref{E:bridge})  that
$D_\infty^{\prime\prime}\tilde \pi_\infty^{(i)}=0$.   In particular, $\tilde\pi_\infty^{(i)}$ defines a saturated  subsheaf $\widetilde E_\infty^{(i)}$ of $E_\infty$.  Furthermore, it is clear that $\rk(\widetilde
E_\infty^{(i)})=\rk(E_\infty^{(i)})$.  Also, we claim that $\deg(\widetilde E_\infty^{(i)})=\deg(E_\infty^{(i)})$.  To see this, note that since 
$D_\infty^{\prime\prime}\tilde \pi_\infty^{(i)}=0$, and $\Lambda F_j\to\Lambda F_\infty$  and $\pi_j^{(i)}\to \tilde\pi_\infty^{(i)}$ in $L^2$:
\begin{align*}
\deg(\widetilde E_\infty^{(i)})&=\frac{1}{2\pi}\int_X \tr\left(\sqrt{-1}\Lambda F_\infty\tilde\pi_\infty^{(i)}\right)\, dvol
=\lim_{j\to \infty}\frac{1}{2\pi}\int_X \tr\left(\sqrt{-1}\Lambda F_j \pi_j^{(i)}\right)\, dvol\\
&=\deg(E_\infty^{(i)})+\frac{1}{2\pi}\lim_{j\to \infty}\Vert D_j^{\prime\prime}\pi_j^{(i)}\Vert^2_{L^2}
=\deg(E_\infty^{(i)})\ ,
\end{align*}
 as claimed.  
 Now the maximal destabilizing subsheaf $\filt_1(E_\infty)$ of $E_\infty$ is the \emph{unique} saturated subsheaf of $E_\infty$ with this rank and slope (cf.\ \cite[Lemma V.7.17]{Ko}).  Hence,
$\tilde\pi_\infty^{(1)}=\pi_\infty^{(1)}$.
Notice also that since $D_j^{\prime\prime}\overset{\scriptscriptstyle L^2_{loc.}}{\lra} D_\infty^{\prime\prime}$, 
  (\ref{E:ffconvergence}) and (\ref{E:bridge}) imply that $\pi_j^{(1)}\overset{\scriptscriptstyle L^2_{1,loc.}}{\lra} \pi_\infty^{(1)}$.  
Proceed by induction as follows:  fix $1\leq k< \ell$, and assume $\tilde \pi_\infty^{(i)}= \pi_\infty^{(i)}
  $ for $i\leq k$.  Then $\widetilde E_\infty^{(k+1)}/E_\infty^{(k)}$ has the same rank and slope as  the maximal destabilizing subsheaf of $E_\infty/E_\infty^{(k)}$, and is
therefore equal to it as above.  Again we conclude that
$\widetilde E_\infty^{(k+1)}=E_\infty^{(k+1)}$.  Continuing until $k=\ell$ completes the proof of part (1) the lemma.

     For part (2), notice that the argument given above applies to a sequence of Seshadri filtrations as well, where because of the lack of uniqueness of Seshadri filtrations we may conclude only that the ranks and degrees of $\widetilde E_\infty^{(i)}$ correspond to those of $ E_\infty^{(i)}$.
\end{proof}

\begin{proof}[Proof of Thm.\ \ref{T:approximate}]  
Let $D_t$ denote a solution to the YM flow equations on $E\to X$ with initial condition $D_0=(\dbar_E,H)$, and let
 $D_\infty$ be the Uhlenbeck limit for some sequence $D_{t_j}$. Then $\Vert F_{D_{t_j}}\Vert_{L^\infty}$ is uniformly bounded, and by Lemma \ref{L:lp}, $\Lambda F_{D_{t_j}}\overset{\scriptscriptstyle L^p}{\lra}\Lambda F_{D_\infty}$  for all $1\leq
p<\infty$. Moreover, 
 we have shown in Thm.\ \ref{T:type} that the HN type of the limit $(E_\infty, D_\infty^{\prime\prime})$ is the same as that of $(E_, D_0^{\prime\prime})$.
 Hence, the \P\ \ref{A:1} (1-3) are satisfied, and we may apply Lemma \ref{L:projectionconvergence} to conclude that
     $\Psi_j\overset{\scriptscriptstyle L^p}{\lra}\Psi_\infty$ for all $p$, where
     $\Psi_j=\Psi_{\omega}^{hn}(D_{t_j}^{\prime\prime}, H)$ and $\Psi_\infty=\Psi_{\omega}^{hn}(D_\infty^{\prime\prime}, H_\infty)$.
      Since $D_\infty$ is Yang-Mills, 
$\sqrt{-1}\Lambda F_{D_\infty}=\Psi_\infty$
(cf.\ Prop.\ \ref{P:split}).  So
$
\Vert \sqrt{-1}\Lambda F_{D_{t_j}}-\Psi_j\Vert_{L^p}\leq \Vert \Lambda F_{D_{t_j}}-\Lambda F_{D_\infty}\Vert_{L^p}
+\Vert \Psi_j-\Psi_\infty\Vert_{L^p}\to 0
$.
\end{proof}

%%%%%%%%%%%%%%%%%%%%%%%%%%%%%%%%%%%%%%%%%%%%%%%%%%%%%%%%%%%%%%%%%%%%%%%%%%%%%%%%%%%%%%%%%%%%%%%%%%%%%%%%%%%%%%%%%%
%%%%%%%%%%%%%%%%%%%%%%%%%%%%%%%%%%%%%%%%%%%%%%%%%%%%%%%%%%%%%%%%%%%%%%%%%%%%%%%%%%%%%%%%%%%%%%%%%%%%%%%%%%%%%%%%%%

\section{Proof of the Main Theorem}   \label{S:main}

%%%%%%%%%%%%%%%%%%%%%%%%%%%%%%%%%%%%%%%%%%%%%%%%%%%%%%%%%%%%%%%%%%%%%%%%%%%%%%%%%%%%%%%%%%%%%%%%%%%%%%%%%%%%%%%%%%
%%%%%%%%%%%%%%%%%%%%%%%%%%%%%%%%%%%%%%%%%%%%%%%%%%%%%%%%%%%%%%%%%%%%%%%%%%%%%%%%%%%%%%%%%%%%%%%%%%%%%%%%%%%%%%%%%%
 In this final section we complete the proof of the Main Theorem.  The missing ingredient is an identification of the holomorphic structure  of the Uhlenbeck limit.  As stated in the Introduction, this part of the argument applies to minimizing sequences as well.  Therefore, both Thm.'s 1 and 2 will follow from:
 
  \begin{Thm}  \label{T:holomorphic}
Let $D_0$ be an integrable unitary connection on $E\to X$, and let $\vec\mu_0$ be the HN type of 
$(E,D_0^{\prime\prime})$.  Suppose that $D_j$ is a sequence of integrable connections 
in the complex gauge orbit of $D_0$ such that:
 $\HYM(D_j)\to \HYM(\vec\mu_0)$.
  Then 
  there is a YM connection $D_\infty$ on a hermitian bundle $E_\infty\to X$ and a finite set of points $\ansing$ such that:
  \begin{enumerate}
  \item  
  $(E_\infty,D_\infty^{\prime\prime})$ is holomorphically isomorphic to $\Grhns(E,D_0^{\prime\prime})^{\ast\ast}$;
  \item  $E$ and $E_\infty$ are  identified outside $\ansing$ via $L^p_{2,loc.}$ isometries for all $p$;
  \item  Via the isometries in (2), and after passing to a subsequence, $D_j\to D_\infty$ in $L^2_{loc.}$ away from $\ansing$.
  \end{enumerate}
\end{Thm}
 
 \noindent
   The main idea for the proof of Thm.\ \ref{T:holomorphic}  follows  Donaldson \cite{Do1} who constructs a nontrivial holomorphic map $(E, D_0^{\prime\prime})\to (E_\infty,D_\infty^{\prime\prime})$.  With such a map in hand, one may then apply 
   the basic principle that a nontrivial holomorphic map between stable bundles of the same rank and degree must be an isomorphism.
   
     First, however, let us  reduce the problem to the case where the Hermitian-Einstein tensors $\Lambda F_{D_j}$ are uniformly bounded.  Let $D_{j,t}$ denote the solution to the YM flow equations with initial conditions $D_j$ at  time $t$.  Fix $t_0>0$.
  It  follows from Lemma \ref{L:subsolution} that 
  $$
 | \Lambda F_{D_{j,t}}|^2(x)\leq \int_X K_t(x,y)  | \Lambda F_{D_{j,t}}|^2(y)\ dvol(y)\ ,
 $$
 where $K_t(x,y)$ is the heat kernel on $X$.  Since $0< K_t(x,y)\leq C(1+t^{-2})$ for some constant $C$ (cf.\ \cite{CL}), it follows that
  for   $t\geq t_0 >0$,
 $\Vert\Lambda F_{D_{j,t}}\Vert_{L^\infty}$ is uniformly bounded for all $j$ in terms of  $\Vert\Lambda F_{D_j}\Vert_{L^2}$.  Note also that 
 $\HYM(\vec\mu_0)\leq \HYM(D_{j,t})\leq \HYM(D_{j,t_0})\leq \HYM(D_j)
$.
Next, fix $ \delta_0>0$.
By Prop.'s \ref{P:key2}, \ref{P:equality}, and Thm.\ \ref{T:type}, it follows that for each $j$
we may find $t_j\geq t_0$ such that:
\begin{equation} \label{E:delta}
\HYM_\alpha(\vec\mu_0)\leq \HYM_\alpha(D_{j,t_j})\leq \HYM_\alpha(\vec\mu_0)+\delta_0\ ,
\end{equation}
for all $1\leq\alpha\leq 2$.  Moreover, as in the proof of Prop.\ \ref{P:convergence}, we may choose the $\{t_j\}$ so  that $\Vert D_{j,t_j}\Lambda F_{D_{j,t_j}}\Vert_{L^2}\to 0$.
 By Prop.\ \ref{P:ymcompactness}, we may assume, after passing to a subsequence, that $D_{j,t_j}$ has an Uhlenbeck limit $D_{\infty}$ which is a Yang-Mills connection on a bundle $E_\infty$, $L^p_{2,loc.}$ isometric to $E$ off a finite set of points $\ansing$.  Moreover, if $\delta_0$ is chosen sufficiently small in (\ref{E:delta}), then the HN type of $(E_\infty, D_{\infty}^{\prime\prime})$ is $\vec\mu_0$ (see Section \ref{S:hartman}).
We now argue as in the proof of Prop.\ \ref{P:uniqueness} (see also (\ref{E:chernclass})):
 $$
2 \Vert D_{j, t_j}-D_j\Vert_{L^2}^2 \leq \HYM(D_j)-\HYM(D_{j,t_j})  
\leq \HYM(D_j)-\HYM(\vec\mu_0)\ .
$$
Since $D_{j,t_j}\overset{\scriptscriptstyle L^p_{loc.}}{\lra} D_\infty$  and $\HYM(D_j)\to\HYM(\vec\mu_0)$,  it follows that $D_j\overset{\scriptscriptstyle L^2_{loc.}}{\lra} D_\infty$.
Therefore, we may assume from the beginning that $\Vert\Lambda F_{D_j}\Vert_{L^\infty}$ is bounded uniformly in $j$.

 Let $\ansing$ denote the bubbling set of Uhlenbeck limit $D_j\weakarrow D_\infty$.
Associated to an initial HNS filtration $\{\pi_0^{(i)}\}$ of $(E, D_0^{\prime\prime})$ there is an algebraic singular set $\algsing$.  Set $Z=\ansing\cup\algsing$ and $\Omega=X\setminus Z$.
Next, we recall from Lemma \ref{L:projectionconvergence} that we may assume there is a sequence $\{\pi_j^{(i)}\}$  of HNS filtrations of $(E, D_j^{\prime\prime})$, with ranks constant in $j$, such that for each $i$, $\pi_j^{(i)}\to \pi_\infty^{(i)}$ in $L^p\cap L^2_{1,loc.}(\Omega)$.  Here, $\{\pi_\infty^{(i)}\}$   is  a filtration of $(E, D_\infty^{\prime\prime})$  on $X$ by holomorphic subbundles with the same ranks and degrees as the $\pi_j^{(i)}$.    We will prove the result inductively on the length of the HNS filtration.  The inductive hypotheses on the bundle $E\to \Omega$ are the following:
\begin{Assume} \label{A3}
\begin{enumerate}
\item $D_j=g_j(D_0)$ on $\Omega$ for complex gauge transformations $g_j$;
\item $D_j\weakarrow D_\infty$ weakly in $L^p_{1,loc.}(\Omega)$, where $D_\infty$ is Yang-Mills;
\item  $(E,D_0^{\prime\prime})$ and $(E_\infty,D_\infty^{\prime\prime})$ extend to $X$  as reflexive sheaves with the same HN type;
\item  $\Lambda F_{D_j}$ is bounded in $L^\infty_{loc.}(\Omega)$ uniformly in $j$.
\end{enumerate}
\end{Assume}
\noindent The conclusion of the inductive argument will be that  $(E_\infty,D_\infty^{\prime\prime})$ 
is holomorphically isomorphic to $\Gr^{hns}(E,D_0^{\prime\prime})^{\ast\ast}$.

To achieve this, let $S\subset (E,D_0^{\prime\prime})$ be the stable subbundle with $\mu(S)=\mu_{max}(E,D_0^{\prime\prime})$ corresponding to the  initial element $\pi_0=\pi_0^{(1)}$ of the filtration $\{\pi_0^{(i)}\}$, and let $Q=E/S$.  It follows from Prop.\ \ref{P:general} (3) that $\Gr^{hns}(E,D_0^{\prime\prime})=S\oplus \Gr^{hns}(Q)$.  
Let $\Omega_0\subset \Omega$ be the complement of a union of balls around the points of $Z$.
By the proof of \cite[Lemma 2.2]{Bu3} (which also works for weak $L^p_1$ convergence, $p> 4$; see also Prop.\ \ref{P:webster}), after passing to a subsequence we may find holomorphic maps $f_j:S\to (E,D_j^{\prime\prime})$ which converge weakly in $L^p_{2, loc.}(\Omega)$ to a nonzero holomorphic map $f_\infty: S\to (E,D_\infty^{\prime\prime})$.  The map $f_\infty$, in turn,  extends to $X$ by Hartog's Theorem.  If $\pi_j$ denotes the projection to  $f_j(S)$, then  as mentioned above $\pi_j\overset{\scriptscriptstyle L^2_{1,loc.}}{\lra}\pi_\infty$, where $\pi_\infty$ is a subbundle of the same rank and degree as $S$, and $\pi_\infty f_\infty=f_\infty$, $D_\infty^{\prime\prime}\pi_\infty=0$.  Write $\Gr^{hns}(E,D_\infty^{\prime\prime})=S_\infty\oplus Q_\infty$, where $S_\infty=\pi_\infty(E)$, $Q_\infty=\ker\pi_\infty$.  
It follows (cf.\ \cite[Cor.\ V.7.12]{Ko}) that $f_\infty$ must be an isomorphism onto its image: $S\to S_\infty\subset E$.  In particular, $f_\infty$ is everywhere injective.  Since $f_j\to f_\infty$ locally uniformly on $\Omega$, it is easy to verify that $\pi_j\to \pi_\infty$ locally  uniformly as well.  By Lemma \ref{L:fundamentalform} we may assume, after passing to a subsequence, that $\pi_j\to\pi_\infty$ weakly in $L^p_{2,loc.}(\Omega)$ and strongly in $L^p_{1,loc.}(\Omega)$ for all $p$.   After applying  a suitable sequence of gauge transformations 
which are uniformly bounded in $L^p_{1,loc.}(\Omega)$ (cf.\ \cite[Lemma 5.12]{D}) we may assume from the beginning that $D_j^{\prime\prime}$ \emph{preserves} the subbundle $S$.  With this understood, we are ready to use induction.  We have shown that the induced connection $\pi_0 D_j\pi_0$ on $S$ converges to 
a connection on $S_\infty$ whose holomorphic structure is isomorphic to $S$.  This is the first step in the induction.  Now consider the induced connections $D_j^Q=\pi_0^\perp D_j \pi_0^\perp$ on $Q$.  These still satisfy the hypotheses \P\ \ref{A3} (1-3) above.  By Lemma \ref{L:fundamentalform}, the second fundamental forms $D_j^{\prime\prime}\pi_0$ are locally uniformly bounded, so by \cite[I.6.12]{Ko}, $D_j^Q$ also satisfy  \P\ \ref{A3} (4).  By induction then,
$Q_\infty\simeq\Gr^{hns}(Q,(D_0^Q)^{\prime\prime})^{\ast\ast}$.  This completes the proof Thm.\ \ref{T:holomorphic}.

\end{document}